\documentclass[a4paper,11pt]{article}
\usepackage{amsmath, amssymb, amsthm}
\usepackage{graphicx} 

\begin{document}

\def\BF{{\mathbb F}} \def\Aut{{\rm Aut}} \def\Ker{{\rm Ker}}
\def\Id{{\rm Id}} \def\II{{\mathcal I}} \def\Inn{{\rm Inn}}
\def\Sym{{\rm Sym}} \def\Im{{\rm Im}} \def\Supp{{\rm Sup}}
\def\Z{{\mathbb Z}} \def\Syl{{\rm Syl}} \def\CC{{\mathcal C}}
\def\D{{\mathbb D}} \def\F{{\mathbb F}} \def\DD{{\mathcal D}}
\def\C{{\mathbb C}} \def\R{{\mathbb R}} \def\S{{\mathbb S}}
\def\PP{{\mathcal P}} \def\Homeo{{\rm Homeo}} \def\MM{{\mathcal M}}
\def\TT{{\mathcal T}} \def\FF{{\mathcal F}} \def\LL{{\mathcal L}}
\def\Out{{\rm Out}}

\title{\bf{Residual $p$ properties of mapping class groups and surface groups}}
 
\author{\textsc{Luis Paris}}

\date{\today}

\maketitle

\begin{abstract}
\noindent
Let $\MM (\Sigma, \PP)$ be the mapping class group of a punctured oriented surface $(\Sigma, \PP)$ 
(where $\PP$ may be empty), and let $\TT_p(\Sigma,\PP)$ be the kernel of the action of $\MM (\Sigma, 
\PP)$ on $H_1 (\Sigma \setminus \PP, \F_p)$. We prove that $\TT_p( \Sigma, \PP)$ is residually $p$. In 
particular, this shows that $\MM (\Sigma, \PP)$ is virtually residually $p$. For a group $G$ we denote 
by $\II_p(G)$ the kernel of the natural action of $\Out (G)$ on $H_1(G, \F_p)$. In order to achieve our 
theorem, we prove that, under certain conditions ($G$ is conjugacy $p$-separable and has Property A), 
the group $\II_p(G)$ is residually $p$. The fact that free groups and surface groups have Property A is 
due to Grossman. The fact that free groups are conjugacy $p$-separable is due to Lyndon and Schupp. 
The fact that surface groups are conjugacy $p$-separable is, from a technical point of view, the 
main result of the paper.
\end{abstract}



\section{Introduction}

A group $G$ is said to be {\it residually nilpotent} if for every $g \in G \setminus \{1\}$ there 
exists a homomorphism $\varphi: G \to N$ such that $N$ is a nilpotent group and $\varphi(g) \neq 1$. 
For $g,h \in G$, we denote by $[g,h] = g^{-1} h^{-1} gh$ the {\it commutator} of $g$ and $h$, and, if 
$K,H$ are two subgroups of $G$, we denote by $[K,H]$ the {\it commutator subgroup} of $K$ and $H$. This 
is the subgroup generated by $\{ [g,h]; g \in K \text{ and } h \in H\}$, and it is normal (resp. 
characteristic) if both, $K$ and $H$, are normal (resp. characteristic). By a {\it filtration} of a 
group $G$ we mean an infinite descending chain
\[
G=G_1 \triangleright G_2 \triangleright G_3 \triangleright \cdots \triangleright G_n \triangleright 
\cdots
\]
of normal subgroups. We say that $\{ G_n\}_{n=1}^{+\infty}$ is a {\it central filtration} if $[G_n, G_m] 
\subset G_{m+n}$ for all $m,n \ge 1$, and we say that $\{ G_n \}_{n=1}^{+\infty}$ is a {\it separating 
filtration} if $\cap_{n=1}^{+\infty} G_n= \{1\}$. It is well-known that a group $G$ is 
residually nilpotent if and only if it admits a separating central filtration. A central filtration $\{ 
G_n \}_{n=1}^{+\infty}$ is called {\it $\Z$-linear} if $G_n/G_{n+1}$ is a free $\Z$-module for all $n \ge 
1$.

\bigskip\noindent
Let $\Sigma_\rho$ be an oriented closed surface of genus $\rho \ge 1$, let $\PP = \{ P_1, \dots, 
P_n\}$ be a (possibly empty) finite collection of punctures in $\Sigma_\rho$, and let $\Homeo (\Sigma_\rho, \PP)$ be 
the group of homeomorphisms $h: \Sigma_\rho \to \Sigma_\rho$ which preserve the orientation and such 
that $h(\PP) = \PP$.  The {\it mapping class group} of $(\Sigma_\rho, \PP)$ is defined to be the group 
$\MM (\Sigma_\rho, \PP) = \pi_0 (\Homeo (\Sigma_\rho, \PP))$ of isotopy classes of elements of $\Homeo 
(\Sigma_\rho, \PP)$. The natural action of $\MM (\Sigma_\rho, \PP)$ on $H_1(\Sigma_\rho \setminus \PP, 
\Z)$ leads to a homomorphism $\tau: \MM (\Sigma_\rho, \PP) \to GL (H_1(\Sigma_\rho \setminus \PP,\Z))$ 
whose kernel is the {\it Torelli group} $\TT (\Sigma_\rho, \PP)$.

\bigskip\noindent
The following result is more or less stated in \cite{Hain1} but, as far as I know, its first complete 
proof is given in \cite{BasLub1} (see also \cite{Farb1}).

\bigskip\noindent
{\bf Theorem 1.1.} {\it The Torelli group $\TT (\Sigma_\rho, \PP)$ has a $\Z$-linear separating central 
filtration. In particular, $\TT (\Sigma_\rho, \PP)$ is residually nilpotent.}

\bigskip\noindent
The filtration of Theorem 1.1 is actually the well-studied Johnson filtration.

\bigskip\noindent
Let $\FF= \{ G_n\}_{n=1}^{+\infty}$ be a central filtration. Then the group commutator induces on
\[
\LL_\FF (G)= \bigoplus_{n=1}^{+\infty} {G_n \over G_{n+1}}
\]
a Lie bracket under which it becomes a graded Lie $\Z$-algebra. If, like in Theorem 1.1, the filtration 
is $\Z$-linear, then $\LL_\FF (G)$ is torsion free and one can apply the Poincar\'e-Birkhoff-Witt 
theorem to embed $\LL_\FF (G)$ into its enveloping algebra.

\bigskip\noindent
Theorem 1.1 also implies that $\TT (\Sigma_\rho, \PP)$ is biorderable (see \cite{Paris1}), which implies 
many other properties (see \cite{MurRhe1}, \cite{DDRW1}).

\bigskip\noindent
Our goal in this paper is to prove a $\F_p$ analogous to Theorem 1.1.

\bigskip\noindent
We say that a group $G$ is {\it residually $p$} if for all $g \in G \setminus \{1\}$ there exists a 
homomorphism $\varphi: G \to P$ such that $P$ is a finite $p$-group and $\varphi(g) \neq 1$. Note that 
residually $p$ groups are residually nilpotent as well as residually finite. A central filtration $\FF= 
\{ G_n \}_{n=1}^{+\infty}$ of a group $G$ is called {\it $\F_p$-linear} if $G_n/G_{n+1}$ is a finite 
dimensional $\F_p$-vector space for all $n \ge 1$. Note that, in that case, $\LL_\FF (G)$ becomes a Lie 
$\F_p$-algebra, and, therefore, one can apply the Poincar\'e-Birkhoff-Witt theorem as well.

\bigskip\noindent
Let $\Sigma_\rho$ be a closed surface of genus $\rho \ge 1$, and let $\PP = \{ P_1, \dots, P_n\}$ be a 
finite collection of punctures in $\Sigma_\rho$. The natural action of $\MM (\Sigma_\rho, \PP)$ on $H_1 
(\Sigma_\rho \setminus \PP, \F_p)$ leads to a homomorphism $\tau_p : \MM (\Sigma_\rho, \PP) \to GL( 
H_1( \Sigma_\rho \setminus \PP, \F_p))$ whose kernel is the {\it congruence subgroup} $\TT_p 
(\Sigma_\rho, \PP)$. This is a finite index subgroup of $\MM (\Sigma_\rho, \PP)$, and is torsion free if 
$p \neq 2$ (see \cite{Ivano3}).

\bigskip\noindent
Our main result is the following.

\bigskip\noindent
{\bf Theorem 1.2.} {\it $\TT_p( \Sigma_\rho, \PP)$ is residually $p$.}

\bigskip\noindent
{\bf Remark.} We will show in Section 2 (see Proposition 2.3) that a finitely generated group which is 
residually $p$ admits a $\F_p$-linear separating central filtration. On the other hand, $\TT_p 
(\Sigma_\rho, \PP)$ is a finite index subgroup of $\MM (\Sigma_\rho, \PP)$ which is finitely generated, 
thus $\TT_p (\Sigma_\rho, \PP)$ is finitely generated. So, by Theorem 1.2, $\TT_p (\Sigma_\rho, \PP)$ 
admits a $\F_p$-linear separating central filtration.

\bigskip\noindent
The proof of Theorem 1.2 is essentially algebraic and uses the following result known as the 
Dehn-Nielsen-Baer theorem (see \cite{Ivano2}, \cite{ZiVoCo1}).

\bigskip\noindent
{\bf Theorem 1.3.} (Dehn, Nielsen \cite{Niels1}, Baer \cite{Baer1}, Magnus \cite{Magnu1}, Zieschang 
\cite{Ziesc1}). {\it Let $\Sigma_\rho$ be a closed oriented surface of genus $\rho \ge 1$, and let $\PP= 
\{ P_1, \dots, P_n \}$ be a finite collection of punctures in $\Sigma_\rho$. Then $\MM (\Sigma_\rho, 
\PP)$ embeds in $\Out (\pi_1 (\Sigma_\rho \setminus \PP))$.
Moreover, if $\PP = \emptyset$, then $\MM( \Sigma_\rho) = \MM(\Sigma_\rho, \emptyset)$ is an
index 2 subgroup of $\Out (\pi_1 (\Sigma_\rho))$.}

\bigskip\noindent
For a group $G$, we denote by $I_p(G)$ the kernel of the natural homomorphism 
$\Aut (G) \to$
\linebreak 
$GL (H_1 (G, 
\F_p))$, and by $\II_p (G)$ the kernel of the natural homomorphism $\Out(G) \to GL (H_1(G, \F_p))$. Note 
that $\Inn (G) \subset I_p(G)$ and $\II_p (G)= I_p(G) / \Inn (G)$. Then Theorem 1.2 is a direct 
consequence of the following.

\bigskip\noindent
{\bf Theorem 1.4.} {\it Suppose $G$ is either a free group or the fundamental group of a closed 
oriented surface. Then $I_p(G)$ and $\II_p (G)$ are both residually $p$.}

\bigskip\noindent
The fact that the outer automorphism group of a free group is virtually residually $p$ has been 
previously proved by Lubotzky \cite{Lubot1}. He did not give explicitely any finite index subgroup 
which is residually $p$, but I suspect he knew that $\II_p (F)$ is an eligible one.
The fact that $I_p(G)$ is residually $p$ is more or less known to experts and follows from the following 
theorem which is proved in Section 2 (see Theorem 2.4).

\bigskip\noindent
{\bf Theorem 1.5.} {\it Let $G$ be a finitely generated and residually $p$ group. Then $I_p(G)$ is 
residually $p$.}

\bigskip\noindent
Now, in order to prove that $\II_p(G)$ is residually $p$ for a given group $G$, we require some restrictions on 
$G$.

\bigskip\noindent
For $g,h \in G$, we use the notation $g \sim h$ to mean that $g$ is conjugate to $h$. We say that a 
group $G$ is {\it conjugacy $p$-separable} if for all $g,h \in G$, either $g \sim h$, or there exists a 
homomorphism $\varphi: G \to P$ such that $P$ is a finite $p$-group and $\varphi(g) \not\sim 
\varphi(h)$. We say that $G$ has {\it Property A} if for every automorphism $\alpha \in \Aut (G)$ such 
that $\alpha(g) \sim g$ for all $g \in G$, we have $\alpha \in \Inn (G)$.

\bigskip\noindent
In Section 2 we prove the following.

\bigskip\noindent
{\bf Theorem 1.6.} {\it Let $G$ be a finitely generated group. If $G$ is conjugacy $p$-separable and has 
Property A, then $\II_p(G)$ is residually $p$.}

\bigskip\noindent
The fact that free groups and fundamental groups of oriented closed surfaces have Property A is 
proved in \cite{Gross1}. The fact that free groups are conjugacy $p$-separable is stated in 
\cite{LynSch1}, Proposition~4.8. However, the proof given in \cite{LynSch1} is imprecise and 
incomplete. We give a complete and detailed proof of this result in Section 3. In Section 4 we prove 
the following.

\bigskip\noindent
{\bf Theorem 1.7.} {\it The fundamental group of a closed oriented surface is conjugacy $p$-separable.}

\bigskip\noindent
From a technical point of view, Theorem 1.7 is the main result of the paper.

\bigskip\noindent
{\bf Acknowledgments.} I started this project with Ben McReynolds. Unfortunately, the distance and the 
activities of both made that I continued alone. I hope he will forgive me, and I thank him for all the 
comments and ideas he gave me. Must of the results of Sections 2 and 3 are derived 
from conversations with Laurent Bartholdi. I thank him for his significant contribution. Finally, I 
thank Mustafa Korkmaz for suggesting me to use topological arguments in the proof of Theorem 4.1, 
and for all the conversations we had on the subject.


\section{Automorphism groups and outer automorphism groups of residually $p$ groups}

We fix a prime number $p\ge 2$. If $K$ and $H$ are two subgroups of a group $G$, we denote by 
$[K,H]^p$ the subgroup generated by $\{ [g,h]; g \in K \text{ and } h \in H\} \cup \{ h^p; h \in H\}$. 
Note that we do not have $[K,H]^p = [H,K]^p$ in general even if $K$ and $H$ are both normal. The {\it 
lower $\F_p$-linear central filtration} of $G$ is defined to be the filtration $\{ \gamma_n^p 
G\}_{n=1}^{+\infty}$ where $\gamma_1^p G=G$ and $\gamma_{n+1}^p G = [G, \gamma_n^p G]^p$ for $n\ge 1$. 
Note that $G/ \gamma_2^p G = H_1(G, \BF_p)$.

\bigskip\noindent
The term ``central'' in the above definition can be used because of the following.

\bigskip\noindent
{\bf Lemma 2.1.} {\it Let $n,m \ge 1$. Then
$[\gamma_m^p G, \gamma_n^p G] \subset \gamma_{n+m}^p G$.}

\bigskip\noindent
{\bf Proof.} We argue by induction on $m$. If $m=1$, then
\[
[\gamma_m^pG, \gamma_n^pG] = [G, \gamma_n^pG] \subset [G, \gamma_n^pG]^p = \gamma_{n+1}^p G\,.
\]
Assume $m>1$. We need to prove that $[g,h] \in \gamma_{m+n}^pG$ for every generator $g$ of 
$\gamma_m^p G$ and every element $h \in \gamma_n^pG$. There are two different cases to consider: (1) 
$g \in [G, \gamma_{m-1}^p G]$; (2) $g=g_1^p$ for some $g_1 \in \gamma_{m-1}^p G$.

\bigskip\noindent
Suppose $g \in [G, \gamma_{m-1}^p G]$. By \cite{MaKaSo1}, Theorem 5.2, we have
\[
[[G,\gamma_{m-1}^pG],\gamma_n^pG] \subset [[\gamma_{m-1}^pG,\gamma_n^pG],G] \cdot [\gamma_n^pG, G], 
\gamma_{m-1}^pG]\,.
\]
By induction we have
\begin{gather*}
[[\gamma_{m-1}^p G, \gamma_n^p G], G] \subset [\gamma_{m+n-1}^p G,G] \subset \gamma_{m+n}^p G\,,\\
[[ \gamma_n^p G, G], \gamma_{m-1}^p G] \subset [\gamma_{n+1}^p G, \gamma_{m-1}^p G] \subset 
\gamma_{m+n}^p G\,,
\end{gather*}
thus
\[
[g,h] \in [[ G, \gamma_{m-1}^p G], \gamma_n^p G] \subset \gamma_{m+n}^p G\,.
\]

\bigskip\noindent
Now, we assume that $g=g_1^p$ for some $g_1 \in \gamma_{m-1}^p G$. Let $g_2$ be another element of 
$\gamma_{m-1}^p G$. We have
\[
[g_1g_2,h] = [g_1,h] \cdot [[ g_1,h], g_2] \cdot [g_2,h]\,,
\]
and, by induction,
\[
[[ g_1,h], g_2] \in [[ \gamma_{m-1}^p G, \gamma_n^p G], \gamma_{m-1}^p G] \subset [ \gamma_{m+n-1}^p G, 
\gamma_{m-1}^p G] \subset \gamma_{2m+n-2}^p G \subset \gamma_{m+n}^p G\,,
\]
thus
\[
[g_1g_2, h] \equiv [g_1,h] \cdot [g_2,h]\text{ mod } \gamma_{m+n}^p G\,.
\]
Iterating this congruence we obtain
\[
[g,h] = [g_1^p,h] \equiv ([g_1,h])^p \text{ mod }\gamma_{m+n}^p G\,,
\]
thus $[g,h] \in \gamma_{m+n}^p G$ since $[g_1,h] \in [\gamma_{m-1}^p G, \gamma_n^pG] \subset 
\gamma_{m+n-1}^p G$.
\qed

\bigskip\noindent
The proof of the following is left to the reader.

\bigskip\noindent
{\bf Lemma 2.2.} {\it Let $P$ be a finitely generated group. Then $P$ is a finite $p$-group if and only 
if there exists some $n \ge 1$ such that $\gamma_n^p P = \{1\}$.}
\qed

\bigskip\noindent
By a {\it $p$-filtration} of a group $G$ we mean a filtration $\{G_n\}_{n=1}^{+\infty}$ such that 
$G/G_n$ is a finite $p$-group for all $n \ge 1$.  
Note that the $p$-filtrations that will appear in the proofs of Theorems 2.4 and 2.5 are actually
$\F_p$-linear filtrations, but this is not needed for our purpose.
Recall also that a filtration $\{ 
G_n\}_{n=1}^{+\infty}$ is said to be a {\it separating filtration} if $\cap_{n=1}^{+\infty} G_n = \{1\}$.

\bigskip\noindent
The main criterion we will use in the proofs of Theorems 2.4 and 2.5 are the following.

\bigskip\noindent
{\bf Proposition 2.3.} {\it
\begin{enumerate}
\item
Let $G$ be a (discrete) group. If $G$ has a separating $p$-filtration, then $G$ is residually $p$.
\item
Let $G$ be a finitely generated group. Then $G$ is residually $p$ if and only if the lower $\F_p$-linear 
central filtration $\{ \gamma_n^p G \}_{n=1}^{+\infty}$ is a separating filtration.
\end{enumerate}}

\bigskip\noindent
{\bf Remark.} It is easily seen that $\gamma_n^pG$ is finitely generated for all $n \ge 1$ and that $\{ 
\gamma_n^p G\}_{n=1}^{+\infty}$ is a $\F_p$-linear filtration if $G$ is finitely generated. 
This is not necessarily true anymore 
if $G$ is not finitely generated.

\bigskip\noindent
{\bf Proof.} Part (1) is easy to proof and left to the reader. So, we assume that $G$ is finitely 
generated and turn to prove Part (2). If $\{ \gamma_n^p G\}_{n=1}^{+\infty}$ is separating, then $G$ is 
residually $p$ because $\{ \gamma_n^p G\}_{n=1}^{+\infty}$ is a $p$-filtration (apply Part (1)).
Suppose that $G$ is residually $p$. Let $g \in G \setminus \{1\}$. There exists a homomorphism $\varphi 
: G \to P$ such that $P$ is a finite $p$-group and $\varphi (g) \neq 1$. By Lemma 2.2, there exists $n 
\ge 1$ such that $\gamma_n^p P= \{1\}$. We have $\varphi( \gamma_n^p G) \subset \gamma_n^p P = \{1\}$, 
thus $g \not\in \gamma_n^p G$ since $\varphi (g) \neq 1$. This shows that $\{ \gamma_n^p 
G\}_{n=1}^{+\infty}$ is a separating filtration.
\qed

\bigskip\noindent
{\bf Theorem 2.4.} {\it Let $G$ be a finitely generated group. If $G$ is residually $p$, then $I_p(G)$ 
is residually $p$.}

\bigskip\noindent
{\bf Proof.} We consider the lower $\F_p$-linear central filtration $\{ \gamma_n^p G\}_{n=1}^{+\infty}$. Since 
$\gamma_{n+1}^p G$ is a characteristic subgroup of $G$, the quotient map $G \to G/ \gamma_{n+1}^p G$ 
determines a homomorphism
\[
\pi_n: \Aut G \to \Aut (G/ \gamma_{n+1}^p G)\,.
\]
We set $A_n= \Ker \pi_n$ for all $n \ge 1$ and turn to prove that $\{A_n\}_{n=1}^{+\infty}$ is a 
separating $p$-filtration of $A_1= I_p(G)$.

\bigskip\noindent
{\bf Claim 1.} {\it $A_n$ acts trivially on $\gamma_k^p G/ \gamma_{n+k}^p G$ for all $k \ge 1$.}

\bigskip\noindent
{\bf Proof of Claim 1.} We argue by induction on $k$. The case $k=1$ follows from the definition of 
$A_n$, so we can assume $k>1$. We take $\alpha \in A_n$, and turn to prove that $\alpha (g) \equiv g 
\text{ mod } \gamma_{n+k}^p G$ for every generator $g$ of $\gamma_k^p G$.

\bigskip\noindent
Suppose first that $g = [g_1, h]$ for some $g_1 \in \gamma_{k-1}^p G$ and $h \in G$. By induction, 
there exist $x_1 \in \gamma_{k+n-1}^p G$ and $y \in \gamma_{n+1}^p G$ such that $\alpha (g_1)= g_1x_1$ 
and $\alpha (h) = hy$. By Lemma 2.1,
\begin{gather*}
[g_1x_1,y] \in [\gamma_{k-1}^pG, \gamma_{n+1}^pG] \subset \gamma_{k+n}^pG\,,\\
[[g_1x_1,h],y] \in [[\gamma_{k-1}^pG,G], \gamma_{n+1}^pG] \subset [\gamma_k^pG, \gamma_{n+1}^pG] 
\subset \gamma_{k+n+1}^pG \subset \gamma_{k+n}^p G\,,
\end{gather*}
thus
\[
\alpha (g) = [g_1x_1, hy] = [g_1x_1, h] \cdot [g_1x_1,y] \cdot [[g_1x_1,h],y] \equiv [g_1x_1,h] \text{ 
mod } \gamma_{k+n}^p G\,.
\]
Again, by Lemma 2.1,
\begin{gather*}
[[g_1,h],x_1] \in [[ \gamma_{k-1}^p G,G] ,\gamma_{k+n-1}^p G] \subset [\gamma_k^p G, \gamma_{k+n-1}^p 
G] \subset \gamma_{2k+n-1}^p G \subset \gamma_{k+n}^pG\,,\\
[x_1,h] \in [\gamma_{k+n-1}^pG,G] \subset \gamma_{k+n}^pG\,,
\end{gather*}
thus
\[
[g_1x_1,h] = [g_1,h] \cdot [[g_1,h],x_1] \cdot [x_1,h] \equiv [g_1,h] = g \text{ mod } \gamma_{k+n}^p 
G\,.
\]
This shows that $\alpha (g) \equiv g \text{ mod } \gamma_{k+n}^pG$.

\bigskip\noindent
Now, we assume that $g= g_1^p$ for some $g_1 \in \gamma_{k-1}^p G$. Let $g_2$ be another element of 
$\gamma_{k-1}^pG$. By induction, there exist $x_1, x_2 \in \gamma_{k+n-1}^p G$ such that $\alpha (g_1) 
= g_1x_1$ and $\alpha(g_2) = g_2x_2$. By Lemma~2.1,
\[
[g_2, x_1^{-1}] \in [\gamma_{k-1}^pG, \gamma_{k+n-1}^pG] \subset \gamma_{2k+n-2}^pG \subset 
\gamma_{k+n}^pG \quad (\text{since } k\ge 2)\,,
\]
thus
\[
\alpha (g_1g_2) = g_1x_1g_2x_2 = g_1g_2[g_2,x_1^{-1}] x_1x_2 \equiv g_1g_2x_1x_2 \text{ mod } 
\gamma_{k+n}^pG\,.
\]
Iterating the above congruence we obtain
\[
\alpha (g) = \alpha (g_1^p) \equiv g_1^p x_1^p \equiv g_1^p \text{ mod } \gamma_{k+n}^pG\,.
\]

\bigskip\noindent
{\bf Claim 2.} {\it The filtration $\{A_n\}_{n=1}^{+\infty}$ is a separating filtration.}

\bigskip\noindent
{\bf Proof of Claim 2.} Let $\alpha \in I_p(G)$, $\alpha \neq \Id$. Choose $g \in G$ such that 
$\alpha(g) \neq g$. Since $\{ \gamma_n^p G\}_{n=1}^{+\infty}$ is a separating filtration, there exists 
$n \ge 1$ such that $g^{-1} \alpha(g) \not\in \gamma_{n+1}^p G$, that is, $\alpha (g) \not\equiv g 
\text{ mod } \gamma_{n+1}^p G$. This shows that $\alpha \not\in A_n$.

\bigskip\noindent
{\bf Claim 3.} {\it Let $\alpha \in A_n$. Then the map
\[\begin{array}{rrcl}
u_\alpha : &G/ \gamma_{n+1}^pG &\to& \gamma_{n+1}^pG / \gamma_{n+2}^pG\\
&[g]&\mapsto&[g^{-1} \alpha(g)]
\end{array}\]
is well-defined.}

\bigskip\noindent
{\bf Proof of Claim 3.} Let $g_1,g_2 \in G$ such that $g_1 \equiv g_2 \text{ mod } \gamma_{n+1}^pG$. 
Let $x_1,x_2 \in \gamma_{n+1}^pG$ such that $\alpha(g_1)= g_1x_1$ and $\alpha(g_2) = g_2x_2$. We have 
to prove that $x_1 \equiv x_2 \text{ mod }\gamma_{n+2}^pG$. Set $a= g_1^{-1}g_2 \in \gamma_{n+1}^pG$. 
By Claim 1, $\alpha(a)$ is of the form $\alpha(a)= ay$ where $y \in \gamma_{2n+1}^pG \subset 
\gamma_{n+2}^pG$. We have
\[
[a,x_1^{-1}] \in [\gamma_{n+1}^pG,\gamma_{n+1}^pG] \subset \gamma_{2n+2}^pG \subset \gamma_{n+2}^pG\,,
\]
and
\[
\alpha(g_2) = g_1x_1ay =g_1 a [a,x_1^{-1}]x_1y = g_2 [a,x_1^{-1}]x_1y\,,
\]
thus
\[
x_2=[a,x_1^{-1}]x_1y \equiv x_1\text{ mod } \gamma_{n+2}^p G\,.
\]

\bigskip\noindent
We denote by $H_n$ the set of set-maps from $G/ \gamma_{n+1}^pG$ to $\gamma_{n+1}^pG/ \gamma_{n+2}^pG$. The 
group $\gamma_{n+1}^pG / \gamma_{n+2}^pG$ is a finite dimensional $\BF_p$-vector space and 
$G/\gamma_{n+1}^pG$ is finite, thus $H_n$ is finite and naturally endowed with a structure of 
$\BF_p$-vector space. The next claim shows that $A_n/A_{n+1}$ is also a finite dimensional $\BF_p$-vector space 
for all $n \ge 1$, so it implies that $\{A_n\}_{n=1}^{+\infty}$ is a $p$-filtration.

\bigskip\noindent
{\bf Claim 4.} {\it the map
\[ \begin{array}{rrcl}
u:&A_n&\to&H_n\\
&\alpha&\mapsto&u_\alpha
\end{array}\]
is a group homomorphism whose kernel is $A_{n+1}$.}

\bigskip\noindent
{\bf Proof of Claim 4.} Let $\alpha_1, \alpha_2 \in A_n$ and let $g \in G$. We write $\alpha_1(g)= 
gx_1$, $\alpha_2(g)= gx_2$, and $\alpha_2(x_1)= x_1y$, where $x_1,x_2\in \gamma_{n+1}^pG$ and $y\in 
\gamma_{2n+1}^pG \subset \gamma_{n+2}^pG$. We have $(\alpha_2\alpha_1) (g)= gx_2x_1y$, thus
\[
u_{\alpha_2 \alpha_1} ([g]) = [x_2x_1y] = [x_2x_1] = u_{\alpha_2} ([g]) \cdot u_{\alpha_1} ([g]) = 
(u_{\alpha_2} \cdot u_{\alpha_1}) ([g])\,.
\]
This shows that $u$ is a homomorphism.

\bigskip\noindent
Let $\alpha \in A_{n+1}$. If $g \in G$, then $\alpha(g) \equiv g \text{ mod } \gamma_{n+2}^pG$, thus 
$u_{\alpha} ([g]) = [1]$. This shows that $u_\alpha$ is trivial, that is, $\alpha$ belongs to the 
kernel of $u$.

\bigskip\noindent
Let $\alpha$ be an element of the kernel of $u$. Let $g \in G$, and let $x= g^{-1} \alpha(g)$. We have 
$u_\alpha ([g]) = [x] = [1]$, thus $x \in \gamma_{n+2}^pG$, hence $\alpha(g) \equiv g \text{ mod } 
\gamma_{n+2}^pG$. This shows that $\alpha \in A_{n+1}$.
\qed

\bigskip\noindent
{\bf Theorem 2.5.} {\it Let $G$ be a finitely generated group. If $G$ is conjugacy $p$-separable and 
has Property A, then $\II_p(G)$ is residually $p$.}

\bigskip\noindent
{\bf Proof.} Recall that $\II_p(G)= I_p(G)/ \Inn (G)$. So, in order to prove Theorem 2.5, it suffices 
to find a $p$-filtration $\{B_n\}_{n=1}^{+\infty}$ of $I_p(G)$ such that $\cap_{n=1}^{+\infty} B_n = 
\Inn (G)$.

\bigskip\noindent
We denote by $\pi_n: I_p(G) \to \Aut (G/ \gamma_{n+1}^pG)$ the restriction of $\pi_n: \Aut(G) \to \Aut (G/ 
\gamma_{n+1}^pG)$ to $I_p(G)$. The kernel of $\pi_n$ is the group $A_n$ defined in the proof of 
Theorem 2.4, and $\{A_n\}_{n=1}^{+\infty}$ is a separating $p$-filtration of $I_p(G)$. Let $S_n$ be the 
set of conjugacy classes in $G/\gamma_{n+1}^pG$, and let $\Sym (S_n)$ denote the symmetric group of 
$S_n$. The group $\Aut(G/\gamma_{n+1}^pG)$ acts naturally on $S_n$, and this action defines a 
homomorphism $\theta_n : \Aut(G/\gamma_{n+1}^pG) \to \Sym (S_n)$. We set $\tilde \pi_n = \theta_n \circ 
\pi_n : I_p(G) \to \Sym (S_n)$, and $B_n= \Ker \tilde \pi_n$.

\bigskip\noindent
{\bf Claim 1.} {\it $\{B_n\}_{n=1}^{+\infty}$ is a $p$-filtration.}

\bigskip\noindent
{\bf Proof of Claim 1.} The fact that $\{A_n\}_{n=1}^{+\infty}$ is a $p$-filtration means that 
$I_p(G)/A_n \simeq \Im \pi_n$ is a finite $p$-group, thus $I_p(G)/B_n \simeq \Im \tilde \pi_n = 
\theta_n( \Im \pi_n)$ is also a finite $p$-group (it is the image of a finite $p$-group under a 
homomorphism).

\bigskip\noindent
{\bf Claim 2.} {\it $\cap_{n=1}^{+\infty} B_n= \Inn (G)$.}

\bigskip\noindent
{\bf Proof of Claim 2.} The inclusion $\Inn (G) \subset \cap_{n=1}^{+\infty} B_n$ is clear, thus we 
only have to prove the reverse inclusion: $\cap_{n=1}^{+\infty} B_n \subset \Inn (G)$. Let $\alpha \in I_p(G)$, 
$\alpha \not \in \Inn (G)$. Since $G$ has Property A, there exists $g \in G$ such that $g$ is not 
conjugate to $\alpha(g)$. Since $G$ is conjugacy $p$-separable, there exists a homomorphism $\varphi: 
G \to P$ such that $P$ is a finite $p$-group, and $\varphi(g) \not\sim \varphi (\alpha (g))$. By Lemma 
2.2, there exists $n \ge 1$ such that $\gamma_{n+1}^pP = \{1\}$. We have $\varphi (\gamma_{n+1}^pG) 
\subset \gamma_{n+1}^pP = \{1\}$, thus $\gamma_{n+1}^pG \subset \Ker \varphi$. Let $\mu_n: G \to G/ 
\gamma_{n+1}^pG$ be the quotient map. The above inclusion and the fact that $\varphi (g) \not\sim 
\varphi( \alpha(g))$ imply that $\mu_n(g) \not\sim \mu_n(\alpha(g))$. This shows that  
$\tilde \pi_n (\alpha)$ does not fix
the conjugacy class of $[g] \in G/ \gamma_{n+1}^pG$, thus $\tilde \pi_n 
(\alpha) \neq \Id$, therefore $\alpha \not\in B_n$.
\qed


\section{Conjugacy $p$-separability of free groups}

We start with an observation that the reader must keep in mind because it will be often used in 
the remainder of the paper.

\bigskip\noindent
{\bf Lemma 3.1.} {\it Let $G$ be a group, and let $K_1, K_2$ be two normal subgroups such that $G/K_1$ 
and $G/K_2$ are finite $p$-groups. Then $G/(K_1 \cap K_2)$ is a finite $p$-group.}

\bigskip\noindent
{\bf Proof.} Let $\varphi_1: G \to G/K_1$ and $\varphi_2: G \to G/K_2$ be the quotient maps. Then $K_1 
\cap K_2$ is the kernel of the homomorphism $\varphi: G \to G/K_1 \times G/K_2$, $g \mapsto 
(\varphi_1(g), \varphi_2(g))$.
\qed

\bigskip\noindent
Now, we prove the following.

\bigskip\noindent
{\bf Theorem 3.2.} {\it Free groups are conjugacy $p$-separable.}

\bigskip\noindent
{\bf Proof.} Let $F=F(X)$ be a free group freely generated by a set $X$. Let $g = 
x_1^{\varepsilon_1} \cdots x_l^{\varepsilon_l}$ be an element of $F$ written in normal form ({\it i.e.} $x_i 
\in X$, $\varepsilon_i \in \{ \pm 1\}$, and $x_{i+1}^{\varepsilon_{i+1}} \neq x_i^{-\varepsilon_i}$ for 
all $i$). The number $l$ is called the {\it word length} of $g$ and is denoted by $|g|_X$, and the set 
$\{x_1, \dots, x_l\} \subset X$ is called the {\it support} of $g$ and is denoted by $\Supp(g)$.

\bigskip\noindent
We take two elements $g,h \in F$ that are not conjugate and we prove by induction on $|g|_X + |h|_X$ 
that there exists a homomorphism $\varphi: F \to P$ such that $P$ is a finite $p$-group and $\varphi 
(g) \not \sim \varphi(h)$.

\bigskip\noindent
{\bf Step 1.} {\it We assume $g=1$.}

\bigskip\noindent
We have $h \neq 1$ and $F$ is residually $p$ (see \cite{Bourb1}, for example), thus there exists a 
homomorphism $\varphi : F \to P$ such that $P$ is a finite $p$-group and $\varphi(h) \neq 1$. We 
obviously have $\varphi (h) \not \sim \varphi(g) = 1$.

\bigskip\noindent
{\bf Conclusion of Step 1.} We can assume $g \neq 1$ and $h \neq 1$.

\bigskip\noindent
{\bf Step 2.} {\it We assume that $|\Supp(g)| = |\Supp(h)| =1$.}

\bigskip\noindent
Then $g$ and $h$ are of the form $g=x^n$ and $h=y^m$ where $x,y \in X$ and $n,m \in \Z \setminus 
\{0\}$. Moreover, $n \neq m$ if $x=y$. Let $q=p^e$ be a power of $p$ strictly greater than $2|n|$ and 
than $2|m|$. Set $P= \Z/q\Z$, and let $\varphi: F \to P$ be the homomorphism defined by $\varphi(x)=1$, 
and $\varphi(z)=0$ for all $z \in X \setminus \{x\}$. We have $\varphi(g) \neq \varphi(h)$, thus 
$\varphi (g) \not \sim \varphi(h)$ since $P$ is abelian.

\bigskip\noindent
{\bf Conclusion of Step 2.} We can assume that $|\Supp (g)| \ge 2$.

\bigskip\noindent
For $f \in F$, we denote by $[f]_1$ the homology class of $f$ in $H_1(F)= \oplus_{x \in X} \Z [x]_1$.

\bigskip\noindent
{\bf Step 3.} {\it We assume that $[g]_1 \neq [h]_1$.}

\bigskip\noindent
Write
\[
[g]_1-[h]_1 = \sum_{i=1}^l a_i [x_i]_1\,,
\]
where $x_1, \dots, x_l \in X$ and $a_i \in \Z \setminus \{0\}$. Let $q=p^e$ be a power of $p$ strictly 
greater than $2|a_1|$. Set $P=\Z/q\Z$, and let $\varphi: F \to P$ be the homomorphism defined by 
$\varphi(x_1) =1$, and $\varphi(z)=0$ for all $z \in X \setminus \{x_1\}$. We have $\varphi(g) \neq 
\varphi(h)$, thus $\varphi(g) \not\sim \varphi(h)$.

\bigskip\noindent
{\bf Conclusion of Step 3.} We can assume that $[g]_1 = [h]_1$.

\bigskip\noindent
{\bf Step 4.} {\it We prove that there exist a homomorphism $\mu: F \to \Z /p\Z$ and an element $x 
\in \Supp (g)$ such that $\mu(x)=1$ and $\mu(g)= \mu(h)=0$.}

\bigskip\noindent
Suppose $[g]_1 = [h]_1 = 0$. We choose some $x \in \Supp(g)$ and we define $\mu: F \to \Z/p\Z$ by 
$\mu(x)=1$, and $\mu(z)=0$ for all $z \in X \setminus \{x\}$. We have $\mu(g)=\mu(h)=0$.

\bigskip\noindent
Suppose $[g]_1 = [h]_1 = a [y]_1$ for some $y \in X$ and $a \in \Z \setminus \{0\}$. Take $x \in 
\Supp(g) \setminus \{y\}$ (this is possible because $|\Supp (g)| \ge 2$) and define $\mu: F \to \Z/p\Z$ 
by $\mu(x)=1$, and $\mu(z)=0$ for all $z \in X \setminus \{x\}$. Then $\mu(g)= \mu(h)= 0$.

\bigskip\noindent
Suppose $[g]_1 = [h]_1 = a_1 [y_1]_1 + \cdots + a_l [y_l]_1$, where $y_1, \dots, y_l \in X$, $a_1, 
\dots, a_l \in \Z \setminus \{0\}$, and $l \ge 2$. Let $m= \gcd (a_1, \dots, a_l)$, $b_i= a_i/m$ for $1 
\le i\le l$, and $T_1= b_1 [y_1]_1 + \dots + b_l [y_l]_1$. Let $U= \Z [y_1]_1 \oplus \cdots \oplus \Z 
[y_l]_1$ be the $\Z$-module freely generated by $[y_1]_1, \dots, [y_l]_1$. Since $\gcd(b_1, \dots,b_l)=1$, 
$T_1$ can be completed into a $\Z$-basis $\{ T_1, T_2, \dots, T_l\}$ of $U$. We define $\bar \mu: 
H_1(F) \to \Z/p\Z$ by
\[
\bar \mu(T_1) = 0\,,\quad \bar\mu (T_2)=1\,, \quad \bar\mu(T_i)=0 \text{ for } 3\le i\le l\,, 
\quad\text{and } \bar\mu(z)=0 \text{ for } z \in X \setminus \{ y_1, \dots, y_l\}\,.
\]
The fact that $\bar\mu(T_1)=0$ implies that $\bar\mu ([g]_1) = \bar\mu ([h]_1) =0$, and the fact that 
$\bar\mu (T_2) \neq 0$ implies that there exists $x \in \{y_1, \dots, y_l\} \subset \Supp (g)$ such 
that $\bar\mu ([x]_1) \neq 0$. Each element of $(\Z/pZ) \setminus \{0\}$ generates $\Z/p\Z$, thus one 
can find an automorphism $\xi: \Z/p\Z \to \Z/p\Z$ such that $\xi (\bar\mu ([x]_1)) =1$. Set $\mu= \xi 
\circ \bar\mu \circ \kappa: F \to \Z/p\Z$, where $\kappa: F \to H_1(F)$ is the quotient map. Then 
$\mu(x)=1$ and $\mu(g)=\mu(h)=0$.

\bigskip\noindent
Now, for $y\in X \setminus \{x\}$ and $0 \le i\le p-1$, we set
\[
d(y,i)= x^iy x^{-\mu(x^iy)}\,.
\]
Then $\Ker\mu$ is the free group $F(Y)$ freely generated by
\[
Y= \{ d(y,i)\ ;\ 0\le i\le p-1\text{ and } y\in X \setminus \{x\} \} \cup \{ x^p\}\,.
\]
For $f \in \Ker \mu = F(Y)$ and $0 \le i \le p-1$, we set $f_i= x^i f x^{-i} \in 
F(Y)$.

\bigskip\noindent
{\bf Step 5.} {\it Let $f \in F(Y)$. We show that $|f_i|_Y \le |f|_X$ for all $0 \le i\le p-1$, and that, if $x \in 
\Supp(f)$, then there exists $i \in \{0,1, \dots, p-1\}$ such that $|f_i|_Y< |f|_X$.}

\bigskip\noindent
Set
\[
d(x,i)= \left\{\begin{array}{ll}
1\quad&\text{if }0\le i\le p-2\\
x^p\quad&\text{if } i=p-1
\end{array}\right.\quad d(x^{-1},i)=\left\{\begin{array}{ll}
x^{-p}\quad&\text{if }i=0\\
1\quad&\text{if } 1\le i \le p-1
\end{array}\right.
\]
For $y \in X \setminus \{x\}$ and $0 \le i\le p-1$ we set
\[
d(y^{-1},i)= x^i y^{-1} x^{-\mu(x^iy^{-1})}\,.
\]
Note that
\[
d(y^{-1},i) = d(y, \mu(x^iy^{-1}))^{-1}
\]
for all $y \in X$ (included $y=x$) and all $0 \le i \le p-1$.

\bigskip\noindent
Let $f= y_1^{\varepsilon_1} \cdots y_l^{\varepsilon_l}$ be the normal form of $f$ with respect to $X$. 
For $1 \le j\le l$ we set
\[
i_j= \mu( x^i y_1^{\varepsilon_1} \cdots y_{j-1}^{\varepsilon_{j-1}})\,.
\]
Then
\[
f_i= d(y_1^{\varepsilon_1},i_1) \cdot d( y_2^{\varepsilon_2},i_2) \cdots d(y_l^{\varepsilon_l},i_l)\,.
\]
This shows that $|f_i|_Y \le |f|_X$.

\bigskip\noindent
Suppose, moreover, that some $y_j$ is equal to $x$. We can choose $i \in \{ 0,1, \dots, p-1\}$ so that 
$i_j \neq p-1$ if $\varepsilon_j=1$, or $i_j \neq 0$ if $\varepsilon_j =-1$. In this case 
$d(y_j^{\varepsilon_j},i_j)=1$ and, therefore, $|f_i|_Y<|f|_X$.

\bigskip\noindent
{\bf Conclusion of Step 5.} We can assume that $|g_{i_0}|_Y<|g|_X$ for some $0\le i_0 \le p-1$, 
and $|h_i|_Y \le |h|_X$ for all $0 \le i \le p-1$.

\bigskip\noindent
By the inductive hypothesis, there exists a homomorphism $\alpha_i: F(Y) \to A_i$ such that $A_i$ is a 
finite $p$-group and $\alpha_i (g_{i_0}) \not\sim \alpha_i (h_i)$, for all $0 \le i\le p-1$. Let $L = 
\cap_{i=0}^{p-1} \Ker \alpha_i$, let $B=F(Y)/L$, and let $\beta: F(Y) \to B$ be the quotient map. Then 
$B$ is a finite $p$-group (by Lemma 3.1) and $\beta(g_{i_0}) \not\sim \beta(h_i)$ for all 
$0 \le i\le p-1$. Let $K= 
\cap_{i=0}^{p-1} x^iLx^{-i}$. It is easily seen that $K$ is a normal subgroup of $F=F(X)$ and that 
$P=F(X)/K$ is a finite $p$-group. We denote by $\varphi: F \to P$ the quotient map.

\bigskip\noindent
{\bf Step 6.} {\it We prove that $\varphi(g) \not\sim \varphi(h)$.}

\bigskip\noindent
Suppose not. Then there exists $f \in F(X)$ such that $\varphi(g_{i_0})= \varphi(fhf^{-1})$. Write $f$ 
in the form $f=f'x^i$ where $f' \in F(Y)$ and $0\le i\le p-1$. Then $\varphi(g_{i_0}) = \varphi(f' h_i 
{f'}^{-1})$ which implies that $\beta(g_{i_0}) = \beta (f'h_i {f'}^{-1})$: a contradiction.
\qed
 

\section{Conjugacy $p$-separability of surface groups}

In this section we prove the following.

\bigskip\noindent
{\bf Theorem 4.1.} {\it The fundamental group of a closed oriented surface of genus $\rho \ge 1$ is 
conjugacy $p$-separable.}

\bigskip\noindent
We fix a genus $\rho \ge 1$ and we denote by $\Sigma_\rho$ the oriented closed surface of genus $\rho$. 
If $\rho = 1$, then $\pi_1 (\Sigma_\rho) \simeq \Z^2$ which is obviously conjugacy $p$-separable, thus 
we can and will assume $\rho \ge 2$.

\bigskip\noindent
We use in our proof the fact that $\pi_1(\Sigma_\rho)$ can be expressed (in many ways) as the 
amalgamated product of two free groups along a cyclic group. So, we start with some preliminaries on 
amalgamated products of groups.

\bigskip\noindent
Let $G= H_1 \ast_C H_2$ be the amalgamated product of two groups $H_1,H_2$ along a subgroup $C$. A 
{\it normal form} of an element $g \in G \setminus C$ is defined to be an expression
\[
g= g_1 \cdot g_2 \cdots g_l
\]
of $g$ such that
\begin{itemize}
\item
either $g_i \in H_1 \setminus C$ or $g_i \in H_2 \setminus C$, for all $1 \le i\le l$,
\item
if $g_i \in H_1 \setminus C$ then $g_{i+1} \in H_2 \setminus C$, and if $g_i \in H_2 \setminus C$ then 
$g_{i+1} \in H_1 \setminus C$, for all $1 \le i\le l-1$.
\end{itemize}
Such an expression always exists but is not unique in general. However, the length $l$ of the normal 
form is unique. It is called the {\it Syllable length} of $g$ and is denoted by $\Syl (g)$. A normal 
form is said to be {\it cyclically reduced} if is length is either $1$ or even. Note that, in the 
latest case (when $l$ is even) we have the equivalences
\[
g_l \in H_2 \setminus C \Leftrightarrow g_1 \in H_1 \setminus C\,, \quad g_l \in H_1 \setminus C 
\Leftrightarrow g_1 \in H_2 \setminus C\,.
\]
If a normal form of $g$ is cyclically reduced, then all its normal forms are cyclically reduced. We say 
in that case that $g$ itself is {\it cyclically reduced}. We also use the convention that the elements 
of $C$ are cyclically reduced and their syllable length are equal to $0$. It is easily checked that 
every element of $G$ is conjugate to a cyclically reduced element.

\bigskip\noindent
We refer to \cite{Serre1} for the proof of the following proposition and for a general exposition on 
amalgamated products of groups.

\bigskip\noindent
{\bf Proposition 4.2.} {\it Let $G = H_1 \ast_C H_2$ be the amalgamated product of two groups $H_1,H_2$ 
along a subgroup $C$, and let $g,h \in G$ be two cyclically reduced elements.
\begin{enumerate}
\item
If $g \sim h$, then $\Syl (g) = \Syl (h)$.
\item
Suppose $\Syl (g) = \Syl (h) = l \ge 2$, and let $h= h_1 \cdots h_l$ be a normal form of $h$. We have 
$g \sim h$ if and only if there exist $u \in C$ and $1 \le t\le l$ such that
\[
ugu^{-1} = h_t \cdot h_{t+1} \cdots h_l \cdot h_1 \cdots h_{t-1}\,.
\]
\item
Suppose $\Syl(g) = \Syl (h)=l \ge 2$, and let $g= g_1 \cdots g_l$ and $h =h_1 \cdots h_l$ be normal 
forms of $g$ and $h$, respectively. There exists $u \in C$ such that $ugu^{-1} =h$ if and only if there 
exist $l+1$ elements $u_0, u_1, \dots, u_l \in C$ such that $u_0=u_l$ and $u_{i-1} g_i = h_i u_i$ for all 
$1 \le i\le l$.
\end{enumerate}}
\qed

\bigskip\noindent
An element $g$ of a group $G$ is said to be {\it conjugacy $p$-distinguished} if, for every $h \in G$, 
either $h \sim g$, or there exists a homomorphism $\varphi: G \to P$ such that $P$ is a finite 
$p$-group and $\varphi(g) \not\sim \varphi(h)$. Note that a group $G$ is conjugacy $p$-separable if and 
only if all its elements are conjugacy $p$-distinguished.

\bigskip\noindent
The strategy we use to prove Theorem 4.1 is the following. We start with an element $g \in \pi_1( 
\Sigma_\rho) \setminus \{1\}$ and show that there exists a decomposition $\pi_1( \Sigma_\rho) = F_1 
\ast_C F_2$ of $\pi_1(\Sigma_\rho)$ as the amalgamated product of two free groups such that $g$ is 
cyclically reduced of syllable length $\ge 2$ (see Proposition 4.4). Afterwards, we take an element $h 
\in \pi_1 (\Sigma_\rho)$ which is not conjugate to $g$, and we prove that there exists a homomorphism 
$\varphi: F_1 \ast_C F_2 \to P_1 \ast_{\bar C} P_2$ such that $P_1 \ast_{\bar C} P_2$ is the amalgamated 
product of two finite $p$-groups along a cyclic subgroup $\bar C$, $\varphi(g)$ is cyclically reduced, 
$\Syl (\varphi(g)) = \Syl(g) \ge 2$, and $\varphi(h) \not \sim \varphi(g)$ (see Proposition 4.5). We 
conclude using the following.

\bigskip\noindent
{\bf Theorem 4.3} (Ivanova \cite{Ivano1}). {\it Let $H= P_1 \ast_{\bar C} P_2$ be the amalgamated 
product of two finite $p$-groups $P_1, P_2$ along a cyclic group $\bar C$. Then every cyclically 
reduced element of $H$ of syllable length $\ge 2$ is conjugacy $p$-distinguished.}
\qed

\bigskip\noindent
{\bf Remark.} The exact result which can be found in \cite{Ivano1} is that: if the amalgamated product 
$P_1 \ast_{\bar C} P_2$ of two finite $p$-groups is residually $p$, then any element of $P_1 \ast_{\bar 
C} P_2$ of infinite order is conjugacy $p$-distinguished. The fact that $g$ is of infinite order if it 
is cyclically reduced and of syllable length $\ge 2$ is a classical fact that can be found for instance 
in \cite{Serre1}. The fact that $P_1 \ast_{\bar C} P_2$ is residually $p$ if $\bar C$ is cyclic is 
proved in \cite{Higma1}. Note that the amalgamated product of two finite $p$-groups is not residually 
$p$ in general (see \cite{Higma1}).

\bigskip\noindent
We turn now to the first step of our proof of Theorem 4.1.

\bigskip\noindent
{\bf Proposition 4.4.} {\it Let $g \in \pi_1( \Sigma_\rho) \setminus \{1\}$. Then there exists a 
decomposition $\pi_1(\Sigma_\rho) = F_1 \ast_C F_2$ of $\pi_1(\Sigma_\rho)$ as an amalgamated product 
such that $F_1 = F(x_1, \dots, x_n,y_1, \dots, y_n)$ is 
a free group of rank $2n$,
$F_2=F(x_1', \dots, x_m', y_1', \dots, y_m')$ is
a free group of rank $2m$, 
$C$ is the cyclic subgroup generated by
\[
\gamma = [x_1,y_1] \cdots [x_n,y_n] = [x_1', y_1'] \cdots [x_m',y_m']\,,
\]
and the element $g \in F_1 \ast_C F_2$ is cyclically reduced of syllable length $\ge2$.}

\bigskip\noindent
{\bf Proof.} Let $a : \S^1 \hookrightarrow \Sigma_\rho$ be a separating simple closed curve of 
$\Sigma_\rho$ (see Figure 4.1). The curve $a$ separates $\Sigma_\rho$ into two subsurfaces $\Sigma_1$ 
and $\Sigma_2$ , each of them with a unique boundary component. Let $a_1: \S^1 \hookrightarrow \partial 
\Sigma_1$ (resp. $a_2: \S^1 \hookrightarrow \partial \Sigma_2$) be the boundary curve of $\Sigma_1$ 
(resp. of $\Sigma_2$). Then $\Sigma_\rho = (\Sigma_1 \sqcup\Sigma_2)/\sim$, where $\sim$ is the 
equivalence relation which identifies $a_1(z)$ with $a_2(\bar z)$ for all $z \in \S^1$. Set $B_0 = 
a(1)$, let $\bar a: [0,1] \to \Sigma_\rho$ be the loop based at $B_0$ defined by $\bar a(t) = a 
(e^{2i\pi t})$, and let $\gamma \in \pi_1 (\Sigma_\rho,B_0)$ be the class of $\bar a$. Then $\pi_1 
(\Sigma_\rho)$ has the amalgamated decomposition $\pi_1 (\Sigma_\rho, B_0) = \pi_1 (\Sigma_1, B_0) 
\ast_C \pi_1( \Sigma_2, B_0)$, where $C$ is the infinite cyclic subgroup generated by $\gamma$. Let $n$ 
(resp. $m$) be the genus of $\Sigma_1$ (resp. $\Sigma_2$). Then $\pi_1 (\Sigma_1, B_0) = F_1 = F(x_1, 
\dots, x_n, y_1, \dots, y_n)$ is a free group of rank $2n$, $\pi_1( \Sigma_2, B_0) = F_2 = F(x_1', 
\dots, x_m', y_1', \dots, y_m')$ is a free group of rank $2m$, and the generators can be chosen so that
\[
\gamma= [x_1, y_1] \cdots [x_n,y_n] = [x_1',y_1'] \cdots [x_m',y_m']\,.
\]

\begin{figure}[htbp]
\bigskip
\centerline{
\setlength{\unitlength}{.5cm}
\begin{picture}(22,6)
\put(0,1){\includegraphics[width=11cm]{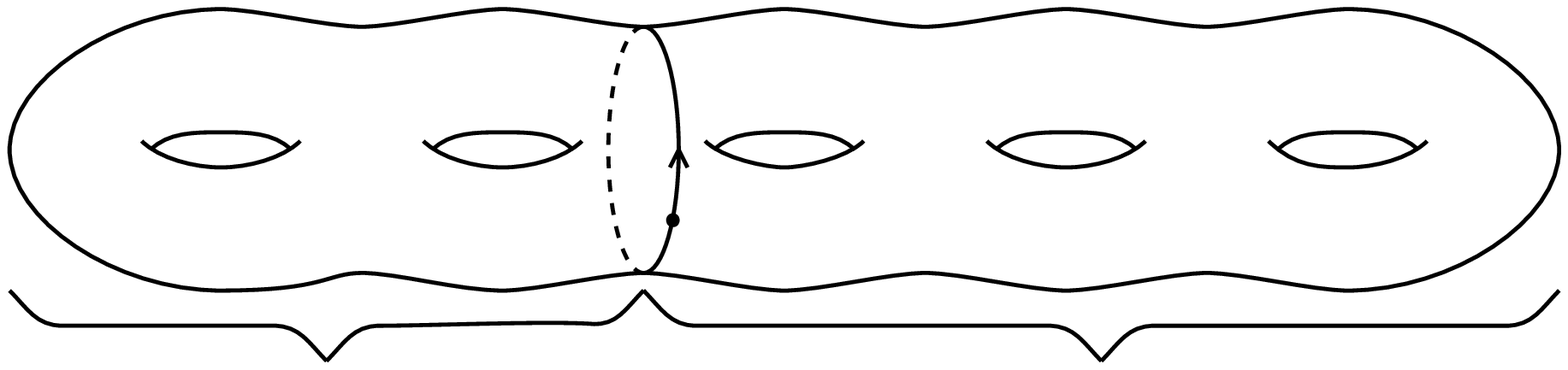}}
\put(4.2,0.2){\small $\Sigma_1$}
\put(15.2,0.2){\small $\Sigma_2$}
\put(9.6,2.8){\small $B_0$}
\put(9.5,5){\small $a$}
\end{picture}}
\bigskip
\centerline{{\bf Figure 4.1.}  Amalgamated decomposition of $\pi_1(\Sigma_\rho)$.}
\end{figure}

\bigskip\noindent
Let $g \in \pi_1( \Sigma_\rho, B_0) \setminus \{1\}$. Let $\bar c_g: [0,1] \to \Sigma_\rho$ be a loop 
based at $B_0$ which represents $g$, and let $c_g: \S^1 \to \Sigma_\rho$ be the closed curve defined 
by $c_g (e^{2i\pi t}) = \bar c_g(t)$. It is easily shown that an element $h \in \pi_1( \Sigma_\rho, 
B_0)$ is conjugate to $g$ if and only if $c_h$ is homotopic to $c_g$ (here the homotopies are free 
with no fixed point).

\bigskip\noindent
Let $c: \S^1 \to \Sigma_\rho$ be a closed curve homotopic to $c_g$ and which intersects $a$ 
transversely. Let $\bar c: [0,1] \to \Sigma_\rho$ be the loop defined by $\bar c(t)= c(e^{2i\pi t})$, 
and let $0 \le t_1<t_2< \cdots <t_l<1$ such that $\bar c(t_i) \in a(\S^1)$ for all $1 \le i\le l$, and 
$\bar c(t) \not\in a(\S^1)$ if $t \not\in \{ t_1, \dots, t_l \}$. In other words, we have $c(\S^1) \cap 
a(\S^1) = \{\bar c(t_1), \dots, \bar c(t_l) \}$, where the points are counted with multiplicity. For $1 
\le i\le l$, let $s_i \in [0,1[$ such that $\bar c(t_i)= \bar a(s_i)$, and let
\[
h_i= \left[ \bar a|_{[0,s_i]} \cdot \bar c|_{[t_i,t_{i+1}]} \cdot \left( \bar a|_{[0,s_{i+1}]} 
\right)^{-1} \right] \in \pi_1 (\Sigma_\rho, B_0)
\]
(here we assume $t_{l+1}=t_1$, $s_{l+1} = s_1$, and $[t_l,t_1]=[t_l,1] \cup [0,t_1]$). We have either $\bar 
c([t_i,t_{i+1}]) \subset \Sigma_1$ or $\bar c([t_i,t_{i+1}]) \subset \Sigma_2$. Moreover, if $\bar 
c([t_i,t_{i+1}]) \subset \Sigma_1$ (resp. $\bar c([t_i,t_{i+1}]) \subset \Sigma_2$), then $h_i \in 
\pi_1(\Sigma_1,B_0)=F_1$ (resp. $h_i \in \pi_1(\Sigma_2,B_0)=F_2$).

\bigskip\noindent
If $h_i \in \langle \gamma \rangle = C$, then there is a closed curve $c'$ homotopic to $c$ such that 
$\bar c'$ follows $\bar c$ between $0$ and $t_i-\varepsilon$, $\bar c'$ follows $\bar c$ between 
$t_{i+1}+\varepsilon$ and $1$, and $\bar c'$ remains in the same component of $\Sigma_\rho \setminus a(\S^1)$ 
between $t_i-\varepsilon$ and $t_{i+1}+\varepsilon$. Hence, we can reduce the number $l$ in this 
manner. So, if $l$ is minimal, then
\[
h=h_1 \cdot h_2 \cdots h_l
\]
is a normal form. It is also cyclically reduced because $l$ is even (since $a$ is separating), and $h$ 
is conjugate to $g$. We set $l=I(a,g)$. By Proposition 4.2, this number is well-defined.

\bigskip\noindent
Now, our goal is to prove that we can choose $a$ so that $I(a,g) \neq 0$. We assume that $I(a,g)=0$ for 
every separating simple closed curve $a: \S^1 \hookrightarrow \Sigma_\rho$, and we look for a 
contradiction.

\bigskip\noindent
 We take $2\rho-3$ separating simple closed curves $a_1, \dots, a_{2\rho-3} : \S^1 \hookrightarrow 
\Sigma_\rho$ that are pairwise disjoint and such that the components of the natural compactification of 
$\Sigma_\rho \setminus (\cup_{i=1}^{2\rho-3} a_i(\S^1))$ are all pantalons (spheres with 3 holes) or 
tori with one hole (see Figure 4.2). Since $I(a_i,g)=0$ for all $1 \le i \le 
2\rho-3$, we can assume that $c$ is included in a component of $\Sigma_\rho 
\setminus (\cup_{i=1}^{2\rho-3} a_i(\S^1))$.

\begin{figure}[htbp]
\bigskip
\centerline{
\setlength{\unitlength}{.5cm}
\begin{picture}(26,8)
\put(0,0){\includegraphics[width=13cm]{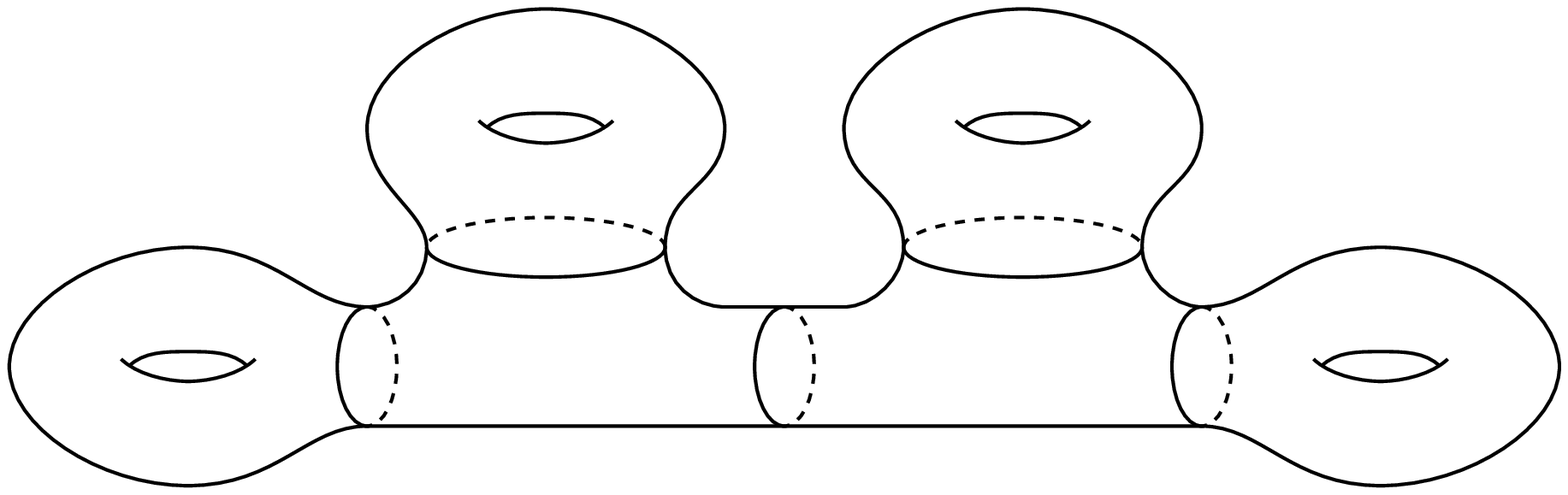}}
\put(4.6,2){\small $a_1$}
\put(8.5,2.8){\small $a_2$}
\put(11.6,2){\small $a_3$}
\put(16.5,2.8){\small $a_4$}
\put(18.6,2){\small $a_5$}
\end{picture}}
\bigskip
\centerline{{\bf Figure 4.2.} A decomposition of $\Sigma_\rho$ into pantalons and one-hole tori.}
\end{figure}

\bigskip\noindent
Let $S_0$ be a tubular neighborhood of $c$, and let $d_1, \dots, d_r : \S^1 \hookrightarrow \partial 
S_0$ be the boundary curves of $S_0$. Up to permutation, we can assume that $d_i$ bounds a disk $\D_i$ 
in $\Sigma_\rho \setminus S_0$ for all $k+1 \le i\le r$, and $d_i$ does not bound any disk for all 
$1 \le i \le k$. Set $S = S_0 \cup (\cup_{i=k+1}^r \D_i)$. Then $S$ is an essential subsurface of 
$\Sigma_\rho$ ({\it i.e.} $\pi_1(S)$ embeds in $\pi_1 (\Sigma_\rho)$) which can be assumed to be 
included in some component of $\Sigma_\rho \setminus (\cup_{i=1}^{2\rho -3} a_i (\S^1))$. In 
particular, $S$ is either an annulus, or a pantalon, or a one-hole torus. 

\bigskip\noindent
We denote by $I(x,y)$ the geometric intersection number of two simple closed curves $x,y$ (see 
\cite{FLP1}). Let $d: \S^1 \hookrightarrow \Sigma_\rho$ be a boundary curve of $S$. Let $a : \S^1 
\hookrightarrow \Sigma_\rho$ be a separating simple closed curve of $\Sigma_\rho$ such that $I(a,d) \ge 
2$ (the proof of the existence of such a curve is left to the reader). Let $a'$ be an arc of $a$ whose 
extremities, $P$ and $Q$, belong to $S$, and let $x$ be a path in $S$ which joins $P$ to $Q$. We 
consider $a' \cup x$ as a closed curve. Then one can prove (with some effort) that: if $a' \cup x$ is 
homotopic to a constant curve, then $a'$ is included in $S$. This fact and the equality $I(a,g)=0$ 
imply that $c$ is homotopic to a curve $c'$ included in $S$ and such that $c' \cap a = \emptyset$. So, 
$c'$ is included in a connected component of $S \setminus a(\S^1)$. Such a component is either a disk 
or an annulus if $S$ is either a pantalon or a one-hole torus, and is a disk if $S$ is an annulus. We 
conclude that $c_g$ is homotopic to a constant curve: a contradiction.
\qed

\bigskip\noindent
The second step of our proof of Theorem 4.1 is given by the following.

\bigskip\noindent
{\bf Proposition 4.5.} {\it Let $\pi_1 (\Sigma_\rho) = F_1 \ast_C F_2$ be an amalgamated decomposition 
of $\pi_1(\Sigma_\rho)$ such that $F_1 = F(x_1, \dots, x_n, y_1, \dots, y_n)$ is a free group of rank 
$2n$, $F_2 = F(x_1', \dots, x_m', y_1', \dots, y_m')$ is a free group of rank $2m$, and $C$ is the 
cyclic group generated by
\[
\gamma = [x_1,y_1] \cdots [x_n,y_n] = [x_1', y_1'] \cdots [x_m', y_m']\,.
\]
Let $g \in \pi_1( \Sigma_\rho)$ be a cyclically reduced element of syllable length $2l \ge 2$, and let 
$h \in \pi_1( \Sigma_\rho)$ be another element not conjugate to $g$. Then there exists a homomorphism 
$\varphi: F_1 \ast_C F_2 \to P_1 \ast_{\bar C} P_2$ such that $P_1 \ast_{\bar C} P_2$ is the 
amalgamated product of two finite $p$-groups along a cyclic group $\bar C$, $\varphi(g)$ is cyclically 
reduced of syllable length equal to $2l= \Syl (g)$, and $\varphi(h)$ is not conjugate to $\varphi(g)$.}

\bigskip\noindent
The following lemmas 4.6 to 4.9 are preliminaries to the proof of Proposition 4.5. Lemmas 4.6 and 4.9 
are due to Stebe (see \cite{Stebe1}, Lemma 1 and Lemma 6). The proof of Lemma 4.8 is long 
and technical, so we put it separately in the last section. 

\bigskip\noindent
{\bf Lemma 4.6} (Stebe \cite{Stebe1}). {\it Let $F=F(X)$ be a free group, let $g \in F \setminus \{1\}$ 
be a non-trivial element, and let $q=p^e$ be a fixed power of $p$. Then there exists a homomorphism 
$\varphi: F \to P$ such that $P$ is a finite $p$-group and the order of $\varphi(g)$ is exactly $q$.}
\qed

\bigskip\noindent
{\bf Lemma 4.7.} {\it Let $F=F(x_1, \dots, x_n, y_1, \dots, y_n)$ be a free group of rank $2n$, let 
$\gamma=$
\linebreak 
$[x_1,y_1] \cdots [x_n,y_n]$, and let $C$ be the cyclic subgroup generated by $\gamma$. Let $g 
\in F \setminus C$. Then there exists a homomorphism $\varphi: F \to P$ such that $P$ is a finite 
$p$-group and $\varphi(g) \not\in \varphi(C)$.}

\bigskip\noindent
{\bf Proof.} Since $\gamma$ is primitive, the centralizer of $\gamma$ in $F$ is equal to $C$. Now, $g 
\not\in C$, thus $[g,\gamma] \neq \{1\}$. The group $F$ is residually $p$, thus there exists a 
homomorphism $\varphi: F \to P$ such that $P$ is a finite $p$-group and $\varphi ([g, \gamma]) = 
[\varphi(g), \varphi(\gamma)] \neq 1$. It follows that $\varphi(g) \not \in \varphi(C) = \langle 
\varphi(\gamma) \rangle$.
\qed

\bigskip\noindent
As pointed out before, the proof of the following is the object of Section 5.

\bigskip\noindent
{\bf Lemma 4.8.} {\it Let $F=F(x_1, \dots, x_n, y_1, \dots, y_n)$ be a free group of rank $2n$, and let 
$\gamma= [x_1, y_1] \cdots [x_n,y_n]$. Let $g,h \in F$. If $\gamma^a g \neq h \gamma^b$ for all $a,b 
\in \Z$, then there exists a homomorphism $\varphi: F \to P$ such that $P$ is a finite $p$-group and 
$\varphi(\gamma)^a \varphi(g) \neq \varphi(h) \varphi(\gamma)^b$ for all $a,b \in \Z$.}

\bigskip\noindent
{\bf Lemma 4.9} (Stebe \cite{Stebe1}). {\it Let $P$ be a finite $p$-group. Let $\xi, \omega \in P$ such that 
$[\omega,\xi \omega \xi^{-1}] \neq 1$, and let $a,b \in \Z$. If $\omega^a = \xi \omega^b \xi^{-1}$, 
then $p$ divides $a$ and $b$.}
\qed

\bigskip\noindent
{\bf Corollary 4.10.} {\it Let $P$ be a finite $p$-group. Let $e \ge 0$, let $\xi, \omega \in P$ such 
that $[\omega^{p^r}, \xi \omega^{p^r} \xi^{-1} ] \neq 1$ for all $0 \le r \le e$, and let $a,b \in \Z$. 
If $\omega^a = \xi \omega^b \xi^{-1}$, then $p^{e+1}$ divides $a$ and $b$.}

\bigskip\noindent
{\bf Proof.} We show by induction on $r$ that $p^{r+1}$ divides $a$ and $b$ for all $0 \le r \le e$. 
The case $r=0$ is a direct consequence of Lemma 4.9. We assume $r \ge 1$. By induction, $p^r$ divides 
$a$ and $b$. We set $a=p^ra'$ and $b=p^rb'$. We have
\[
(\omega^{p^r})^{a'} = \xi (\omega^{p^r})^{b'} \xi^{-1}\,,
\]
thus, by Lemma 4.9, $p$ divides $a'$ and $b'$, therefore $p^{r+1}$ divides $a$ and $b$.
\qed

\bigskip\noindent
{\bf Proof of Proposition 4.5.} We can assume that
\[
g=g_1 \cdot g_1' \cdots g_l \cdot g_l'
\]
is a normal form of $g$, where $g_i \in F_1 \setminus C$ and $g_i' \in F_2 \setminus C$ for all $1 \le 
i\le l$.

\bigskip\noindent
{\bf Step 1.} {\it We assume that $h \in F_1$.}

\bigskip\noindent
By Lemma 4.7, there exists a homomorphism $\alpha_{1\,i} : F_1 \to A_{1\,i}$ such that $A_{1\,i}$ is a 
finite $p$-group and $\alpha_{1\,i} (g_i) \not\in \alpha_{1\,i} (C)$ for all $1 \le i\le l$. Let $L_1 = 
\cap_{i=1}^l \Ker \alpha_{1\,i}$, let $B_1= F_1/L_1$, and let $\beta_1: F_1 \to B_1$ be the quotient 
map. Then $B_1$ is a finite $p$-group (see Lemma 3.1) and $\beta_1(g_i) \not\in \beta_1(C)$ for all $1 \le i\le l$. 
Similarly, there exists a homomorphism $\beta_2: F_2 \to B_2$ such that $B_2$ is a finite $p$-group and 
$\beta_2(g_i') \not\in \beta_2(C)$ for all $1 \le i\le l$. Let $q_1 = p^{e_1}$ be the order of 
$\beta_1(\gamma)$, let $q_2 = p^{e_2}$ be the order of $\beta_2(\gamma)$, and let $q=q_1q_2 = 
p^{e_1+e_2}$. By Lemma 4.6, there exists a homomorphism $\psi_j: F_j \to Q_j$ such that $Q_j$ is a 
finite $p$-group and the order of $\psi_j(\gamma)$ is equal to $q$, for $j=1,2$. Let $K_j= \Ker \beta_j 
\cap \Ker \psi_j$, let $P_j = F_j/K_j$, and let $\varphi_j : F_j \to P_j$ be the quotient map, for 
$j=1,2$. Then $P_1$ is a finite $p$-group, $\varphi_1(g_i) \not\in \varphi_1(C)$ for all $1 \le i\le l$, 
and the order of $\varphi_1(\gamma)$ is equal to $q$. Similarly, $P_2$ is a finite $p$-group, 
$\varphi_2(g_i') \not\in \varphi_2(C)$ for all $1 \le i\le l$, and the order of $\varphi_2 (\gamma)$ is 
equal to $q$. Let $H= P_1 \ast_{\varphi_1(\gamma)= \varphi_2(\gamma)} P_2 = P_1 \ast_{\bar C} P_2$, 
where $\bar C$ is the cyclic group generated by $\varphi_1( \gamma)= \varphi_2 (\gamma)$. Then the 
homomorphisms $\varphi_1, \varphi_2$ induce a homomorphism $\varphi: F_1 \ast_C F_2 \to P_1 \ast_{\bar 
C} P_2$.

\bigskip\noindent
By construction,
\[
\varphi(g) = \varphi_1(g_1) \cdot \varphi_2(g_1') \cdots \varphi_1(g_l) \cdot \varphi_2(g_l')
\]
is a normal form which is cyclically reduced. Furthermore, we have $\varphi(h) \in P_1$, thus, by 
Proposition 4.2, $\varphi(h)$ is not conjugate to $\varphi(g)$.

\bigskip\noindent
{\bf Conclusion of Step 1.} We can assume that $h$ is not conjugate to an element of $F_1$, and is not 
conjugate to an element of $F_2$.
So, we can assume that $h$ is cyclically reduced of syllable length $2k \ge 2$, and has 
a normal form
\[
h=h_1 \cdot h_1' \cdots h_k \cdot h_k'
\]
with $h_j \in F_1 \setminus C$ and $h_j' \in F_2 \setminus C$ for all $1 \le j\le k$.

\bigskip\noindent
{\bf Step 2.} {\it We assume that $k \neq l$.}

\bigskip\noindent
It is easily proved using the same arguments as in Step 1 that there exist homomorphisms $\varphi_1: F_1 
\to P_1$, $\varphi_2 : F_2 \to P_2$ such that $P_1$ and $P_2$ are finite $p$-groups, $\varphi_1(g_i), 
\varphi_1(h_j) \not\in \varphi_1(C)$ for all $1 \le i\le l$ and all $1 \le j\le k$, $\varphi_2 (g_i'), 
\varphi_2( h_j') \not\in \varphi_2(C)$ for all $1 \le i\le l$ and all $1 \le j\le k$, and the order of 
$\varphi_1(\gamma)$ is equal to the order of $\varphi_2(\gamma)$. Let $H= P_1 \ast_{\varphi_1 (\gamma) 
= \varphi_2(\gamma)} P_2 = P_1 \ast_{\bar C} P_2$, where $\bar C$ is the cyclic group generated by 
$\varphi_1 (\gamma) = \varphi_2 (\gamma)$. Then $\varphi_1$ and $\varphi_2$ induce a homomorphism 
$\varphi: F_1 \ast_C F_2 \to P_1 \ast_{\bar C} P_2$.

\bigskip\noindent
By construction
\[
\varphi(g) = \varphi_1 (g_1) \cdot \varphi_2 (g_1') \cdots \varphi_1 (g_l) \cdot \varphi_2(g_l') \quad 
\text{and} \quad \varphi(h) = \varphi_1 (h_1) \cdot \varphi_2 (h_1')\cdots \varphi_1 (h_k) \cdot 
\varphi_2(h_k')
\]
are cyclically reduced normal forms. By Proposition 4.2, we conclude that $\varphi(g) \not\sim 
\varphi(h)$.

\bigskip\noindent
{\bf Conclusion of Step 2.} We can assume $k=l$.

\bigskip\noindent
For the next step we need the following.

\bigskip\noindent
{\bf Claim 1.} {\it Let $f \in F_1 \setminus C$ and $n \in \Z \setminus \{0\}$. Then $[ \gamma^n, f 
\gamma^n f^{-1}] \neq 1$.}

\bigskip\noindent
{\bf Proof of Claim 1.} Suppose that $[\gamma^n, f \gamma^n f^{-1}] =1$. Let $H$ be the subgroup of 
$F_1$ generated by $\{ \gamma, f\}$. Either $H$ is free of rank 2 freely generated by $\{\gamma, f\}$, 
or $H$ is cyclic. Since $[ \gamma^n, f\gamma^n f^{-1}] = 1$, the group $H$ is not freely generated by 
$\{ \gamma, f\}$, thus $H$ is a cyclic group. The element $\gamma$ is primitive, thus $H$ is equal to 
$\langle \gamma \rangle = C$, therefore $f \in C$: a contradiction.

\bigskip\noindent
{\bf Step 3.} {\it We construct a homomorphism $\varphi : F_1 \ast_C F_2 \to P_1 \ast_{\bar C} P_2$, 
where $P_1,P_2$ are finite $p$-groups, $\bar C$ is the cyclic group generated by $\varphi(\gamma)$, and 
$\Syl (\varphi(g))= \Syl(\varphi(h)) = 2l$.}

\bigskip\noindent
For $1 \le t\le l$, we set
\[
h^{(t)} = h_{t+1} \cdot h_{t+1}' \cdots h_{t+l} \cdot h_{t+l}'\,,
\]
where the indices are considered mod $l$. By Proposition 4.2, the fact that $g$ is not conjugate to $h$ 
is equivalent to
\[
\gamma^a h^{(t)} \gamma^{-a} \neq g \quad\text{for all }a \in \Z \text{ and all } 1 \le t\le l\,.
\]
It is easily proved using the same arguments as in Step 1 that there exist homomorphisms $\alpha_1: F_1 
\to A_1$, $\alpha_2: F_2 \to A_2$ such that $A_1,A_2$ are finite $p$-groups, $\alpha_1(g_i), \alpha_1( 
h_i) \not\in \alpha_1 (C)$ for all $1 \le i\le l$, and $\alpha_2 (g_i'), \alpha_2 (h_i') \not\in 
\alpha_2(C)$ for all $1 \le i\le l$.

\bigskip\noindent
Now, for a given $1\le t\le l$, we construct homomorphisms $\beta_1^{(t)} : F_1 \to B_1^{(t)}$, 
$\beta_2^{(t)} : F_2 \to B_2^{(t)}$, where $B_1^{(t)}, B_2^{(t)}$ are finite $p$-groups.

\bigskip\noindent
{\bf Case 1.} There exists $1 \le i\le l$ such that $\gamma^a h_{t+i} \neq g_i \gamma^b$ for all $a,b 
\in \Z$.

\bigskip\noindent
By Lemma 4.8, there exists a homomorphism $\beta_1^{(t)} : F_1 \to B_1^{(t)}$ such that $B_1^{(t)}$ is 
a finite $p$-group and $\beta_1^{(t)} (\gamma^a h_{t+i}) \neq \beta_1^{(t)} (g_i \gamma^b)$ for all 
$a,b \in \Z$. On the other hand, we set $B_2^{(t)} = \{1\}$, and $\beta_2^{(t)}: F_2 \to B_2^{(t)}$ is 
the trivial map.

\bigskip\noindent
{\bf Case 2.} Not in Case 1, and there exists $1 \le i\le l$ such that $\gamma^a h_{t+i}' \neq g_i' 
\gamma^b$ for all $a,b \in \Z$.

\bigskip\noindent
By Lemma 4.8, there exists a homomorphism $\beta_2^{(t)}: F_2 \to B_2^{(t)}$ such that 
$B_2^{(t)}$ is a finite $p$-group and $\beta_2^{(t)} 
(\gamma^a h_{t+i}') \neq \beta_2^{(t)} (g_i' \gamma^b)$ for all $a,b \in \Z$. On the other hand, we set 
$B_1^{(t)} = \{1\}$, and $\beta_1^{(t)}: F_1 \to B_1^{(t)}$ is the trivial map.

\bigskip\noindent
{\bf Case 3.} For all $1 \le i\le l$ there exist $a_i, a_i', b_i, b_i' \in \Z$ such that $\gamma^{a_i} 
h_{t+i} = g_i \gamma^{b_i}$ and $\gamma^{a_i'} h_{t+i}' = g_i' \gamma^{b_i'}$.

\bigskip\noindent
Set $c_i= b_i - a_i'$ for $1 \le i\le l$, $c_i' = b_i' -a_{i+1}$ for $1 \le i\le l-1$, and $c_l'= b_l'-
a_1$. Then
\[
\gamma^{a_1} h^{(t)} \gamma^{-a_1} = g_1 \gamma^{c_1} g_1' \gamma^{c_1'} \cdots g_l \gamma^{c_l} g_l' 
\gamma^{c_l'}\,.
\]
Moreover, since $h$ is not conjugate to $g$, not all the $c_i$'s or the $c_i'$'s are zero. Let
\[
u= \gcd (c_1, \dots, c_l, c_1', \dots, c_l')\,,
\]
and write $u=p^ev$, where $p$ does not divide $v$. Since $F_1$ is residually $p$, by Claim 1, there 
exists a homomorphism $\beta_1^{(t)}: F_1 \to B_1^{(t)}$ such that $B_1^{(t)}$ is a finite $p$-group 
and
\[
\beta_1^{(t)} ([ \gamma^{p^r}, g_i^{-1} \gamma^{p^r} g_i]) \neq 1
\]
for all $1 \le i\le l$ and all $0 \le r\le e$. Similarly, there exists a homomorphism $\beta_2^{(t)} : 
F_2 \to B_2^{(t)}$ such that $B_2^{(t)}$ is a finite $p$-group and
\[
\beta_2^{(t)} ([ \gamma^{p^r}, {g_i'}^{-1} \gamma^{p^r} g_i' ]) \neq 1
\]
for all $1 \le i\le l$ and all $0 \le r\le e$.

\bigskip\noindent
Now, let $L_j= \Ker \alpha_j \cap (\cap_{t=1}^l \Ker \beta_j^{(t)})$, let $N_j= F_j/L_j$, and let $\psi_j : 
F_j \to N_j$ be the quotient map, for $j=1,2$. Let $q_1 = p^{e_1}$ be the order of $\psi_1(\gamma)$, 
let $q_2=p^{e_2}$ be the order of $\psi_2(\gamma)$, and let $q= q_1q_2 = p^{e_1+e_2}$. By Lemma 4.6, 
there exist homomorphisms $\psi_1': F_1 \to N_1'$, $\psi_2': F_2 \to N_2'$ such that $N_1'$ and $N_2'$ 
are finite $p$-groups, and the order of $\psi_j'(\gamma)$ is equal to $q$ for $j=1,2$.

\bigskip\noindent
Let $K_j= \Ker \psi_j \cap \Ker \psi_j'$, let $P_j = F_j/K_j$, and let $\varphi_j: F_j \to P_j$ be the 
quotient map, for $j=1,2$. The groups $P_1,P_2$ are finite $p$-groups, and the orders of 
$\varphi_1(\gamma)$ and $\varphi_2(\gamma)$ are both equal to $q= p^{e_1+e_2}$. Let $H= P_1 
\ast_{\varphi_1(\gamma) = \varphi_2(\gamma)} P_2 = P_1 \ast_{\bar C} P_2$, where $\bar C$ is the cyclic 
group generated by $\varphi_1(\gamma) = \varphi_2 (\gamma)$. Then $\varphi_1$ and $\varphi_2$ induce a 
homomorphism $\varphi: F_1 \ast_C F_2 \to P_1 \ast_{\bar C} P_2$.

\bigskip\noindent
{\bf Step 4.} {\it We prove that $\varphi(h) \not\sim \varphi(g)$.}

\bigskip\noindent
By construction,
\[ \varphi(g)= \varphi_1(g_1)\cdot \varphi_2(g_1') \cdots \varphi_1(g_l) \cdot \varphi_2(g_l') \quad 
\text{and} \quad \varphi(h)= \varphi_1(h_1) \cdot \varphi_2(h_1') \cdots \varphi_1(h_l) \cdot 
\varphi_2(h_l')
\]
are cyclically reduced normal forms. So, by Proposition 4.2, it suffices to show that
\[
\varphi(g) \neq \varphi(\gamma^a h^{(t)} \gamma^{-a})
\]
for all $a \in \Z$ and all $1 \le t\le l$.

\bigskip\noindent
{\bf Case 1.} There exists $1 \le i \le l$ such that $\gamma^a h_{t+i} \neq g_i \gamma^b$ for all $a,b 
\in \Z$.

\bigskip\noindent
By construction of $\beta_1^{(t)}$, we have
\[
\varphi_1 (\gamma^a h_{t+i}) \neq \varphi_1 (g_i \gamma^b)
\]
for all $a,b \in \Z$. By Proposition 4.2, it follows that $\varphi (\gamma^a h^{(t)} \gamma^{-a}) \neq 
\varphi(g)$ for all $a \in \Z$.

\bigskip\noindent
{\bf Case 2.} Not in Case 1, and there exists $1 \le i\le l$ such that $\gamma^a h_{t+i}' \neq g_i' 
\gamma^b$ for all $a,b \in \Z$.

\bigskip\noindent
Then it is easily proved using the same arguments as in Case 1 that $\varphi( \gamma^a h^{(t)} 
\gamma^{-a}) \neq \varphi(g)$ for all $a \in \Z$.

\bigskip\noindent
{\bf Case 3.} For all $1 \le i\le l$ there exist $a_i, a_i', b_i, b_i' \in \Z$ such that $\gamma^{a_i} 
h_{t+i} = g_i \gamma^{b_i}$ and $\gamma^{a_i'} h_{t+i}' = g_i' \gamma^{b_i'}$.

\bigskip\noindent
Following the construction of the homomorphisms $\beta_1^{(t)}$ and $\beta_2^{(t)}$ we set
\[
\gamma^{a_1} h^{(t)} \gamma^{-a_1} = g_1 \gamma^{c_1} g_1' \gamma^{c_1'} \cdots g_l \gamma^{c_l} g_l' 
\gamma^{c_l'}\,,
\]
and $u=p^ev = \gcd (c_1, \dots, c_l, c_1', \dots, c_l')$, where $p$ does not divide $v$. Suppose there 
exists $d \in \Z$ such that $\varphi (\gamma^d h^{(t)} \gamma^{-d}) = \varphi(g)$. By Proposition 4.2, 
there exist $d_1, \dots, d_l, d_1', \dots, d_l' \in \Z$ such that
\[
\varphi (\gamma^{d_{i-1}'} g_i \gamma^{c_i}) = \varphi (g_i \gamma^{d_i}) \quad \text{and}\quad \varphi 
(\gamma^{d_i} g_i' \gamma^{c_i'}) = \varphi(g_i' \gamma^{d_i'})
\]
for all $1 \le i\le l$, where $d_0' = d_l'$. Let $1 \le i\le l$. The above equality implies that
\[
\varphi_1( g_i^{-1} \gamma^{d_{i-1}'} g_i) = \varphi_1 (\gamma^{d_i-c_i})\,.
\]
Moreover, by construction of $\beta_1^{(t)}$,
$\varphi_1( [\gamma^{p^r}, g_i^{-1} \gamma^{p^r} g_i]) \neq 1$
for all $0 \le r \le e$,
thus, by Corollary~4.10, $p^{e+1}$ divides $d_{i-1}'$ and $d_i-c_i$.
Similarly, $p^{e+1}$ divides $d_i$ and $d_i'-c_i'$ for all 
$1 \le i\le l$. It follows that $p^{e+1}$ divides all the $c_i$'s and all the $c_i'$'s: a 
contradiction. So, $\varphi (\gamma^d h^{(t)} \gamma^{-d}) \neq \varphi(g)$ for all $d \in \Z$.
\qed

\bigskip\noindent
{\bf Proof of Theorem 4.1.} Let $g,h \in \pi_1 (\Sigma_\rho)$ such that $g \not\sim h$. We can 
assume $g \neq 1$. By Proposition 4.4, the group $\pi_1 (\Sigma_\rho)$ 
has an amalgamated decomposition $\pi_1 (\Sigma_\rho) = F_1 \ast_C F_2$, where $F_1 = F(x_1, \dots, 
x_n, y_1, \dots, y_n)$ is a free group of rank $2n$, $F_2=F(x_1', \dots, x_m', y_1', \dots, y_m')$ is a 
free group of rank $2m$, and $C$ is the cyclic group generated by
\[
\gamma= [x_1,y_1] \cdots [x_n,y_n] = [x_1',y_1'] \cdots [x_m',y_m']\,,
\]
and such that $g$ is cyclically reduced of syllable length $\ge 2$. By Proposition 4.5, there exists a 
homomorphism $\alpha: F_1 \ast_C F_2 \to P_1 \ast_{\bar C} P_2$ such that $P_1,P_2$ are finite 
$p$-groups, $\bar C$ is the cyclic group generated by $\alpha(\gamma)$, $\alpha(g)$ is cyclically reduced of 
syllable length $\ge 2$, and $\alpha(g) \not\sim \alpha(h)$. By Theorem 4.3, there exists a 
homomorphism $\beta: P_1 \ast_{\bar C} P_2 \to P$ such that $P$ is a finite $p$-group and $\beta( 
\alpha(g)) \not\sim \beta(\alpha (h))$. Let $\varphi= \beta \circ \alpha : \pi_1 (\Sigma_\rho) \to P$. 
Then $\varphi(g) \not\sim \varphi(h)$.
\qed


\section{A technical lemma}

Let $F = F(x_1, \dots, x_n, y_1, \dots, y_n)$ a free group of rank $2n$, and let $\gamma= [x_1, y_1] 
\cdots [x_n,y_n]$. In this section we prove the following.

\bigskip\noindent
{\bf Lemma 4.8.} {\it Let $g,h \in F$. If $\gamma^a g \neq h \gamma^b$ for all $a,b \in \Z$, then there 
exists a homomorphism $\varphi: F \to P$ such that $P$ is a finite $p$-group and $\varphi (\gamma^a g) 
\neq \varphi( h \gamma^b)$ for all $a,b \in \Z$.}

\bigskip\noindent
The following lemmas 5.1 and 5.2 are preliminaries to the proof of Lemma 4.8.

\bigskip\noindent
{\bf Lemma 5.1.} {\it Let $g \in F$. Then there exist a homomorphism $\varphi: F \to \Z /p^2\Z$, and a 
free generating set $Z$ for $\Ker \varphi$, such that $\gamma \in Z$, and $g \gamma g^{-1}$ is of the form 
$g \gamma g^{-1} = w z w^{-1}$, with $w \in \Ker \varphi$ and $z \in Z$.}

\bigskip\noindent
{\bf Proof.} Let $\Sigma$ be a surface of genus $n$ with one hole, and let $c: \S^1 \hookrightarrow 
\partial \Sigma$ be the boundary curve of $\Sigma$. Let $B_0 = c(1)$. Then we identify $F$ with $\pi_1 
(\Sigma, B_0)$ and $\gamma$ with the homotopy class of the loop $\bar c: [0,1] \to \Sigma$, $ t \mapsto 
c(e^{2i\pi t})$.

\bigskip\noindent
Let $\tilde \Sigma$ be the surface of genus $p^2(n-1)+1$ with $p^2$ holes which we represent as follows.

\bigskip\noindent
Choose a small $\varepsilon >0$. For $1 \le k \le n-1$ and $0 \le l \le p^2-1$ we set
\[
O_{k\,l} = k e^{2i\pi l/p^2}\,, \quad V_{k\,l} = (k-1+\varepsilon) e^{2i\pi l/p^2}\,, \quad U_{k\,l} = 
(k-\varepsilon) e^{2i\pi l/p^2}\,.
\]
These are complex numbers. We denote by $\CC_{k\,l}$ the circle of radius $\varepsilon$ centered at 
$O_{k\,l}$, and by $\DD_{k\,l}$ the interval $\DD_{k\,l} = [V_{k\,l}, U_{k\,l}]$, for $1 \le k\le 
n-1$ and $0 \le l \le p^2-1$. We denote by $\CC_0$ the circle centered at $0$ of radius 
$\varepsilon$. We consider the connected graph
\[
\Gamma= \CC_0 \cup (\cup_{k,l} \CC_{k\,l}) \cup (\cup_{k,l} \DD_{k\,l})
\]
embedded in $\C$ (see Figure 5.1).

\begin{figure}[htbp]
\bigskip
\centerline{
\setlength{\unitlength}{.5cm}
\begin{picture}(9,9)
\put(0,0){\includegraphics[width=4.5cm]{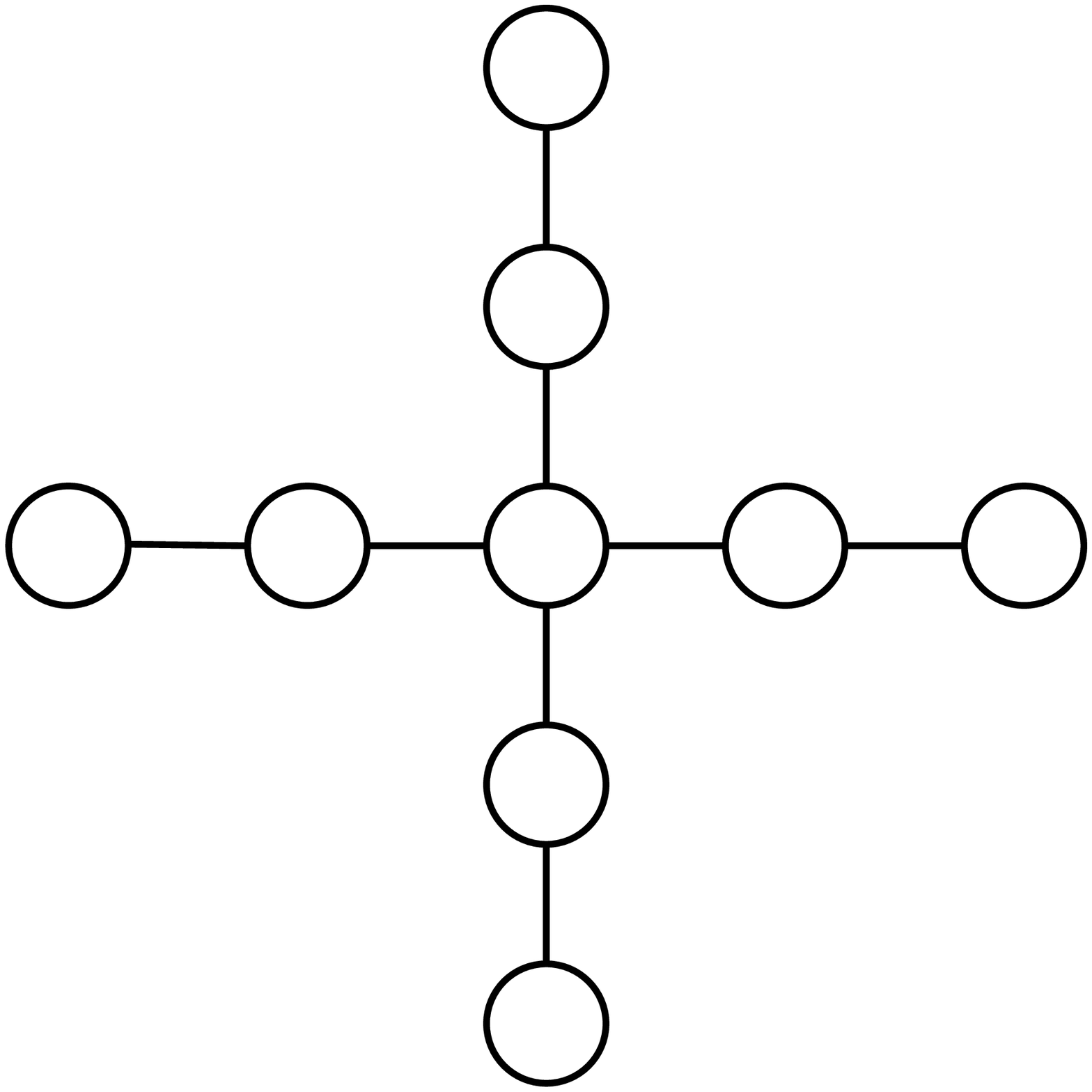}}
\end{picture}}
\bigskip
\centerline{{\bf Figure 5.1.} The graph $\Gamma$ for $n=3$ and $p=2$.}
\end{figure}

\bigskip\noindent
We embed $\C$ in $\R^3= \C \times \R$ by $\xi \mapsto (\xi,0)$, and we denote by $\tilde \Sigma_0$ the 
boundary of a regular neighborhood of $\Gamma$ in $\R^3$. Then $\tilde \Sigma_0$ is a closed oriented 
surface of genus $p^2(n-1)+1$. Let $H_l$ be the plane orthogonal to $\R e^{2i\pi l/p^2}$ and which 
contains $(n-1+\varepsilon -\varepsilon^2) e^{2i\pi l/p^2}$. It is easily seen that this construction 
can be made so that $H_l \cap \tilde \Sigma_0$ is a circle, $\tilde c_l$, which bounds a disk $\D_l$ 
embedded in $\tilde \Sigma_0$. We set $\tilde \Sigma = \tilde \Sigma_0 \setminus (\cup_{l=0}^{p^2-1} 
\D_l)$. Then $\tilde \Sigma$ is an oriented surface of genus $p^2(n-1)+1$ with $p^2$ holes, and its 
boundary curves are $\tilde c_0, \tilde c_1, \dots, \tilde c_{p^2-1}$.

\bigskip\noindent
Let $D= \{0\} \times \R$ be the central vertical line, and let $\theta: \R^3 \to \R^3$ be the rotation 
of angle $2\pi/p^2$ around $D$. Clearly, we can assume that $\theta( \tilde \Sigma) = \tilde \Sigma$, 
and $\theta \circ \tilde c_l = \tilde c_{l+1}$ for all $0 \le l \le p^2-1$ (where the indices are considered 
mod $p^2$). Set $G= \langle \theta \rangle \simeq \Z/ p^2\Z$. By the above, $G$ acts freely on $\tilde 
\Sigma$ and $\tilde \Sigma /G$ is a one-hole surface of genus $n$ which can be identified with 
$\Sigma$. Set $\tilde B_l= \tilde c_l(1)$ for all $0 \le l \le p^2-1$. So, we have a regular covering 
$\pi: \tilde \Sigma \to \Sigma$ (see Figure 5.2) which gives rise to an exact sequence
\[
1 \to \pi_1(\tilde \Sigma, \tilde B_0) \longrightarrow \pi_1 (\Sigma, B_0) \longrightarrow G \simeq 
\Z/p^2\Z \to 1 \,.
\]

\begin{figure}[htb]
\bigskip
\centerline{
\setlength{\unitlength}{.5cm}
\begin{picture}(19,20.5)
\put(0,0){\includegraphics[width=9.5cm]{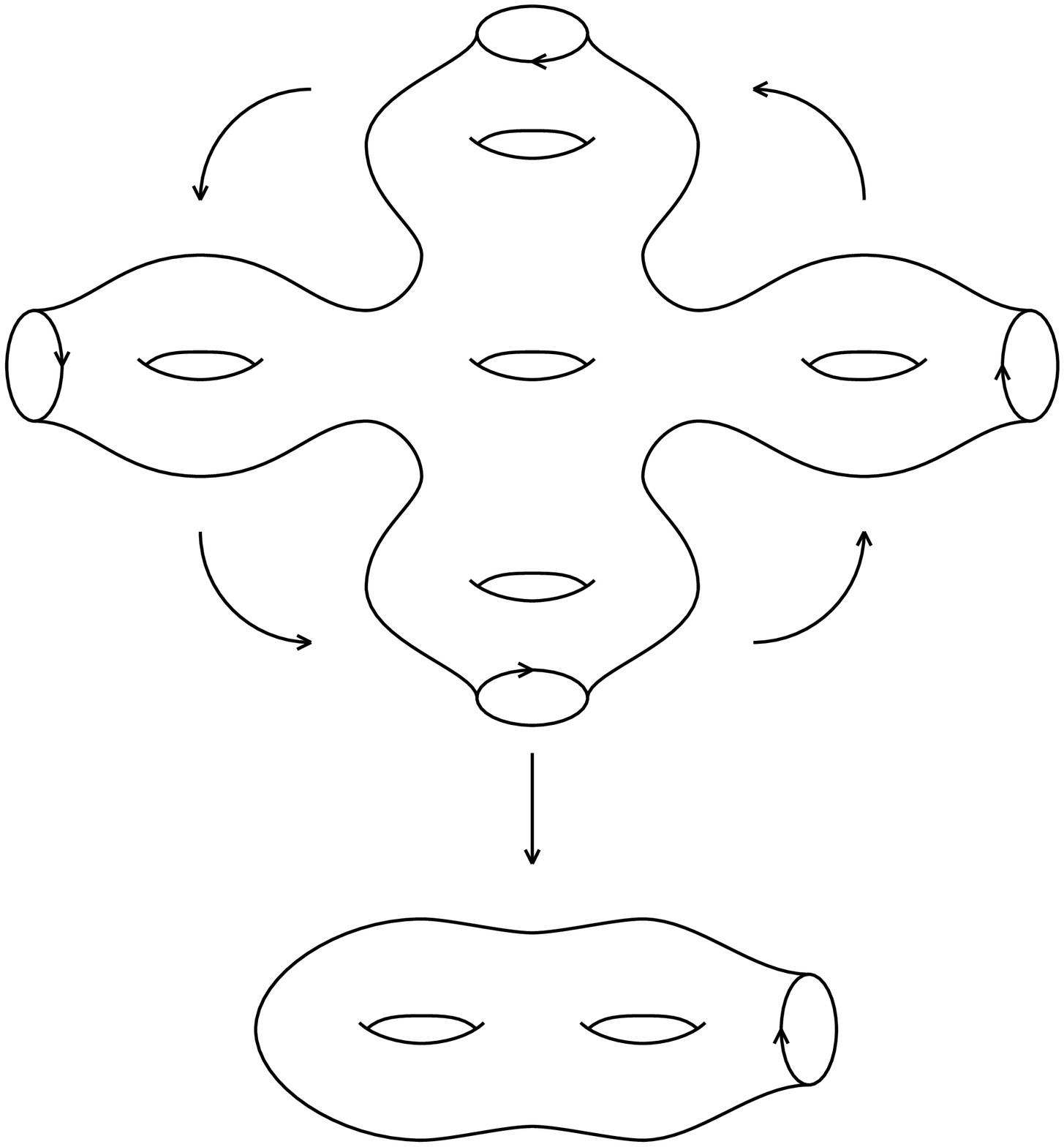}}
\put(13.5,2.2){\small $c$}
\put(9.8,6){$\pi$}
\put(10,8.5){\small $\tilde c_0$}
\put(1.3,13.5){\small $\tilde c_3$}
\put(17.3,14.3){\small $\tilde c_1$}
\put(8.5,19){\small $\tilde c_2$}
\end{picture}}
\bigskip
\centerline{{\bf Figure 5.2.} The regular covering $\tilde \Sigma \to \Sigma$ (for $n=2$ and $p=2$).}
\end{figure}

\bigskip\noindent
For $0 \le l \le p^2-1$, we denote by $\hat c_l : [0,1] \to \tilde \Sigma$ the loop based at $\tilde 
B_l$ defined by $\hat c_l(t)= \tilde c_l (e^{2i\pi t})$, we choose a path $\hat d_l : [0,1] \to \tilde 
\Sigma$ from $\tilde B_0$ to $\tilde B_l$, and we set $z_l = [ \hat d_l \hat c_l \hat d_l^{-1}] \in 
\pi_1( \tilde \Sigma, \tilde B_0)$. For a good choice of the $\hat d_l$'s, the set $\{ z_0, z_1, \dots, 
z_{p^2-1}\}$ can be completed into a generating system $\hat Z= \{ u_1, \dots, u_{p^2(n-1)+1}, v_1, 
\dots, v_{p^2(n-1)+1}, z_0, z_1, \dots, z_{p^2-1} \}$ for $\pi_1(\tilde \Sigma, \tilde B_0)$
so that $\pi_1(\tilde \Sigma, \tilde B_0)$ has a 
presentation with generating set $\hat Z$ and unique relation
\[
[u_1,v_1] \cdots [u_{p^2(n-1)+1}, v_{p^2(n-1)+1}] = z_0z_1 \cdots z_{p^2-1}\,.
\]
In particular, for a chosen $k_0 \in \{ 0,1, \dots, p^2-1\}$, $\pi_1(\tilde \Sigma, \tilde B_0)$ is the 
free group freely generated by $Z = \hat Z \setminus \{z_{k_0}\}$. Also, the path $\hat d_0$ can and 
will be assumed to be the constant path.

\bigskip\noindent
Let $\omega: [0,1] \to \Sigma$ be a loop based at $B_0$ which represents $g$, let $\hat \omega : [0,1] 
\to \tilde \Sigma$ be the lift of $\omega$ such that $\hat \omega (0) = \tilde B_0$, and let $l_0 \in 
\{0,1, \dots, p^2-1\}$ such that $\hat \omega (1)= \tilde B_{l_0}$. Since $p^2 \ge 3$, we can choose 
$k_0$ such that $k_0 \neq 0$ and $k_0 \neq l_0$. Then $\gamma= [c] = [\hat c_0] = z_0 \in Z$, and
\[
g \gamma g^{-1} = [\hat \omega \hat c_{l_0} \hat \omega^{-1}] = w z_{l_0} w^{-1}\,,
\]
where $w= [\hat \omega \hat d_{l_0}^{-1}] \in \pi_1(\tilde \Sigma, \tilde B_0)$.
\qed

\bigskip\noindent
{\bf Remark.} In the above proof we only need the condition $p^2 \ge 3$, so we can replace $p^2$ by $p$ 
if $p \neq 2$.

\bigskip\noindent
{\bf Lemma 5.2.} {\it Let $g,h \in F$. If $\gamma^a g \gamma^b g^{-1} \neq h$ for all $a,b \in \Z$, 
then there exists a homomorphism $\varphi: F \to P$ such that $P$ is a finite $p$-group and $\varphi( 
\gamma^a g \gamma^b g^{-1}) \neq \varphi (h)$ for all $a,b \in \Z$.}

\bigskip\noindent
{\bf Proof.}
By Lemma 5.1, there exist a homomorphism $\pi: F \to \Z/p^2\Z$ and a free generating set $Z$ for 
$\Ker\pi$ such that $\gamma \in Z$ and $g \gamma g^{-1}$ is of the form $u z_0 u^{-1}$, where $u \in 
\Ker\pi$ and $z_0 \in Z$.

\bigskip\noindent
Suppose $\pi(h)\neq 1$. Set $P= \Z/p^2\Z$ and $\varphi = \pi: F \to P$. Then $\varphi( \gamma^a g 
\gamma^b g^{-1}) =1 \neq \varphi(h)$ for all $a,b \in \Z$. So, we can assume that $h \in \Ker\pi$.

\bigskip\noindent
Let $u= \omega_1^{c_1} \omega_2^{c_2} \cdots \omega_l^{c_l}$ be the normal form of $u$ with respect to 
$Z$, where $c_i \in \Z \setminus \{0\}$, $\omega_i \in Z$, and $\omega_i \neq \omega_{i+1}$ for all 
$i$. For convenience, we assume that $\omega_1 = \gamma$ and that $c_1$ may be equal to $0$. Moreover, 
we can also assume that $\omega_l \neq z_0$. We set $v= \omega_2^{c_2} \cdots \omega_l^{c_l}$. For 
$a,b \in \Z$, the normal form of $\gamma^{a+c_1} v z_0^b v^{-1}$ is
\[
\begin{array}{ll}
\gamma^{a+c_1}\quad&\text{if }b=0\,,\\
\gamma^{a+b+c_1}\quad&\text{if } v=1 \text{ and }z_0=\gamma\,,\\
\gamma^{a+c_1} \omega_2^{c_2} \cdots \omega_l^{c_l} z_0^b \omega_l^{-c_l} \cdots \omega_2^{-c_2} 
\quad&\text{otherwise}\,.
\end{array}
\]
Let $h \gamma^{c_1} = z_1^{d_1} \cdots z_k^{d_k}$ be the normal form of $h \gamma^{c_1}$ with respect to 
$Z$. Let $q=p^e$ be a power of $p$ strictly greater than $2|c_i|$ and than $2|d_j|$ for all $2 \le i\le 
l$ and all $1 \le j\le k$. 

\bigskip\noindent
Let
\[
G= \langle Z \ |\ z^q = 1 \text{ for }z \in Z\rangle\,.
\]
The group $G$ is the free product of $|Z|$ copies of $\Z/q\Z$, and there is a natural homomorphism 
$\alpha: \Ker \pi = F(Z) \to G$ which sends $z$ to $z$ for all $z \in Z$. Let
\[
U= \{ \gamma^{a+c_1} v z_0^b v^{-1}; a,b \in \Z\}\,.
\]
Note that the condition $\gamma^a g \gamma^b g^{-1} \neq h$ for all $a,b \in \Z$ is equivalent to $h 
\gamma^{c_1} \not\in U$. Now, observe that $\alpha (h \gamma^{c_1}) = z_1^{b_1} \cdots z_k^{b_k}$ is the 
normal form of $\alpha( h \gamma^{c_1})$ which, by construction, does not coincide with the normal form 
of any element of $\alpha(U)$, thus $\alpha (h \gamma^{c_1}) \not\in \alpha (U)$.

\bigskip\noindent
The set $\alpha (U)$ is finite and the group $G$ is residually $p$, thus there exists a homomorphism 
$\beta : G \to B$ such that $B$ is a finite $p$-group and $(\beta \circ \alpha) (h \gamma^{c_1}) 
\not\in (\beta \circ \alpha) (U)$. Set $\psi= \beta \circ \alpha: \Ker \pi = F(Z) \to B$. Then $\psi(h) 
\neq \psi (\gamma^a u z_0^b u^{-1})$ for all $a,b \in \Z$.

\bigskip\noindent
Let $f \in F$ such that $\pi(f) =1$. Let $K= \cap_{j=0}^{p^2-1} (f^j \Ker \psi f^{-j})$. 
The group $K$ is a normal subgroup of $F$ and the quotient $P=F/K$ is a finite $p$-group. Let $\varphi: 
F \to P$ be the quotient map. Then $\varphi( \gamma^a g \gamma^b g^{-1}) \neq \varphi (h)$ for all $a,b 
\in \Z$.
\qed

\bigskip\noindent
{\bf Proof of Lemma 4.8.} The condition $\gamma^a g \neq h \gamma^b$ for all $a,b \in \Z$ is equivalent 
to $\gamma^a g \gamma^{-b} g^{-1} \neq hg^{-1}$ for all $a,b \in \Z$. By Lemma 5.2, there exists a 
homomorphism $\varphi: F \to P$ such that $P$ is a finite $p$-group and $\varphi (\gamma^a g \gamma^{-
b} g^{-1}) \neq \varphi( hg^{-1})$ for all $a,b \in \Z$. Then $\varphi (\gamma^a g) \neq \varphi( h 
\gamma^b)$ for all $a,b \in \Z$.
\qed



\bigskip\bigskip\noindent
{\bf Luis Paris},

\smallskip\noindent 
Institut de Math\'ematiques de Bourgogne, UMR 5584 du CNRS, Universit\'e de Bourgogne, B.P. 
47870, 21078 Dijon cedex, France

\smallskip\noindent
E-mail: {\tt lparis@u-bourgogne.fr}


\begin{thebibliography}{99}

\bibitem{FLP1}
{\it Travaux de Thurston sur les surfaces.}
S\'eminaire Orsay.
Ast\'erisque, 66-67.
Soci\'et\'e Math\'ematique de France, Paris, 1979.

\bibitem{Baer1}
{\bf R. Baer.}
{\it Isotopie von Kurven auf orientierbaren, geschlossenen Fl\"achen und ihr Zusamen hang mit
der topologischen Deformation der Fl\"achen.}
J. Reine Angew. Math. {\bf 159} (1928), 101-116.

\bibitem{BasLub1}
{\bf H. Bass, A. Lubotzky.}
{\it Linear-central filtrations on groups.}
The mathematical legacy of Wilhelm Magnus: groups, geometry and special functions (Brooklyn, NY, 1992),
45--98, Contemp. Math., 169, Amer. Math. Soc., Providence, RI, 1994.

\bibitem{Bourb1}
{\bf N. Bourbaki.}
{\it Groupes et alg\`ebres de Lie. Chapitres II et III.}
Hermann, Paris, 1972.

\bibitem{DDRW1}
{\bf P. Dehornoy, I. Dynnikov, D. Rolfsen, B. Wiest.}
{\it Why are braids orderable?}
Panoramas et Synth\`eses, 14.
Soci\'et\'e Math\'ematique de France, Paris, 2002.

\bibitem{Farb1}
{\bf B. Farb.}
{\it Some problems on mapping class groups and moduli space.}
Problems on mapping class groups and related topics, 11-55, Proc. Sympos. Pure Math., 74, Amer. Math. 
Soc., Providence, RI, 2006.

\bibitem{Gross1}
{\bf E.K. Grossman.}
{\it On the residual finiteness of certain mapping class groups.}
J. London Math. Soc. (2) {\bf 9} (1974/75), 160--164.

\bibitem{Hain1}
{\bf R. Hain.}
{\it Infinitesimal presentations of the Torelli groups.}
J. Amer. Math. Soc. {\bf 10} (1997), no. 3, 597--651.

\bibitem{Higma1}
{\bf G. Higman.}
{\it Amalgams of $p$-groups.}
J. Algebra {\bf 1} (1964), 301--305.

\bibitem{Ivano3}
{\bf N.V. Ivanov.}
{\it Subgroups of Teichm\"uller modular groups.}
Translations of Mathematical Monographs, 115.
American Mathematical Society, Providence, RI, 1992.

\bibitem{Ivano2}
{\bf N.V. Ivanov.}
{\it Mapping class groups.}
Handbook of geometric topology, 523--633, North-Holland, Amsterdam, 2002.

\bibitem{Ivano1}
{\bf E.A. Ivanova.}
{\it On the approximability with respect to conjugacy of free products of two groups with an
amalgamated subgroup by finite $p$-groups.} (Russian)
Mat. Zametki {\bf 76} (2004), no. 4, 502--509.
Translation in Math. Notes {\bf 76} (2004), no. 3-4, 465--471.

\bibitem{Lubot1}
{\bf A. Lubotzky.}
{\it Normal automorphisms of free groups.}
J. Algebra {\bf 63} (1980), no. 2, 494--498.

\bibitem{LynSch1}
{\bf R.C. Lyndon, P.E. Schupp.}
{\it Combinatorial group theory.}
Reprint of the 1977 edition.
Classics in Mathematics.
Springer-Verlag, Berlin, 2001.

\bibitem{Magnu1}
{\bf W. Magnus.}
{\it \"Uber Automorphismen von Fundamentalgruppen berandeter Fl\"achen.}
Math. Ann. {\bf 109} (1934), 617--646.

\bibitem{MaKaSo1}
{\bf W. Magnus, A. Karrass, D. Solitar.}
{\it Combinatorial group theory.
Presentations of groups in terms of generators and relations.}
Second revised edition.
Dover Publications, Inc., New York, 1976.

\bibitem{MurRhe1}
{\bf R. Botto Mura, A. Rhemtulla.}
{\it Orderable groups.}
Lecture Notes in Pure and Applied Mathematics, Vol. 27.
Marcel Dekker, Inc., New York-Basel, 1977.

\bibitem{Niels1}
{\bf J. Nielsen.}
{\it Untersuchungen zur Topologie der geschlossenen zweiseitigen Fl\"achen.}
Acta Math. {\bf 50} (1927), 189--358.

\bibitem{Paris1}
{\bf L. Paris.}
{\it On the fundamental group of the complement of a complex hyperplane arrangement.}
Arrangements - Tokyo 1998, 257--272, Adv. Stud. Pure Math., 27, Kinokuniya, Tokyo, 2000.

\bibitem{Serre1}
{\bf J.-P. Serre.}
{\it Arbres, amalgames, $SL_2$.}
Ast\'erisque, No. 46.
Soci\'et\'e Math\'ematique de France, Paris, 1977.

\bibitem{Stebe1}
{\bf P.F. Stebe.}
{\it Conjugacy separability of certain free products with amalgamation.}
Trans. Amer. Math. Soc. {\bf 156} (1971), 119--129.

\bibitem{Ziesc1}
{\bf H. Zieschang.}
{\it Discrete groups of plane motions and plane group images.} (Russian)
Uspehi Mat. Nauk {\bf 21} (1966), no. 3, 195--212.

\bibitem{ZiVoCo1}
{\bf H. Zieschang, E. Vogt, H.-D. Coldewey.}
{\it Surfaces and planar discontinuous groups.}
Lecture Notes in Mathematics, 835.
Springer, Berlin, 1980.

\end{thebibliography}
\end{document}